\documentclass[12pt]{article}
    \usepackage{amsmath}
    \usepackage{amssymb}
    \usepackage{amsthm}

\usepackage{pstricks}

      \newcommand {\al}   {\alpha}          \newcommand {\bt}  {\beta}
      \newcommand {\gam } {\gamma}          
      \newcommand {\del}  {\delta}          \newcommand {\Del} {\Delta}
              \newcommand {\ve}   {\varepsilon}
                 \newcommand {\vphi} {\varphi}
      \newcommand {\lam}  {\lambda}         \newcommand {\Lam}  {\Lambda}
      
                \newcommand {\Om}  {\Omega}
      \newcommand {\pl}   {\partial}        \newcommand {\s}    {\sigma}

      \newcommand {\RRR}  {{\mathbb R}}     \newcommand {\SSS}  {{\mathcal S}}
              \newcommand {\MMM}  {{\cal M}}
           \newcommand {\QQQ}  {Q}
      \newcommand {\EEE}  {{\cal E}}        \newcommand {\PPP}  {{\cal P}}
      \newcommand {\FFF}  {{\cal F}}        \newcommand {\BBB}  {{\cal B}}
      \newcommand {\TTT}  {{\mathbb T}}     \newcommand {\NNN}  {\mathcal{N}}
      \newcommand {\Arg} {\text{Arg}}      \newcommand {\interval}{[-\pi/2,\, \pi/2]}
     \newcommand {\beq}  {\begin{equation}}
      \newcommand {\eeq}  {\end{equation}}  \newcommand {\lab}  {\label}
      \newcommand {\set}  {(\pl B \times S^1)_+}
     \newcommand {\iset}  {(I \times S^1)_+}
    \newcommand {\duga} {\stackrel{\smile}}
    \newcommand {\marka} {\text{'}}

      \newtheorem{lem}{Lemma}
      \newtheorem{utv}{Proposition}
      
      \newtheorem{opr}{Definition}
      \newtheorem{predl}{Statement}

\author{Alexander Plakhov}

\title{Billiard scattering on rough sets:\\ Two-dimensional case}

\date{Institute of Mathematical and Physical Sciences,\\
Aberystwyth University, Aberystwyth SY23 3BZ, UK\thanks{On leave
from Department of Mathematics, Aveiro University, Aveiro 3810,
Portugal}}

\begin{document}

\maketitle

\begin{abstract}
The notion of a rough two-dimensional (convex) body is introduced,
and to each rough body there is assigned a measure on $\TTT^3$
describing  billiard scattering on the body. The main result is
characterization of the set of measures generated by rough bodies.
This result can be used to solve various problems of least
aerodynamical resistance.
\end{abstract}

\begin{quote}
\end{quote}

\begin{quote}
\end{quote}

\begin{quote}
{\small {\bf Running title:} Scattering on rough sets}
\end{quote}

\section{Definition of a rough set and statement\\ of main theorem}

\subsection{Introductory remarks and review of literature}

In this paper the notion of a rough two-dimensional (convex) body is
given and some properties of rough bodies are established.

Let $B \subset \RRR^2$ be a convex bounded set with nonempty
interior, that is, a bounded convex body. Consider the "set"{}
obtained from $B$ by moving off a set of "very small"{} area. Such a
(heuristically defined) set is called a {\it rough body}: from the
"macroscopic"{} point of view, it almost coincides with $B$, and
from the "microscopic"{} point of view, it contains some "flaws".
(One can imagine a detail of a mechanism that, after a period of
exploitation, has got some defects.) If the removed set adjoins the
boundary $\pl B$, one can expect that a flow of point particles
incident on the rough body is reflected in another way as compared
to reflection from $B$.

The notion of rough body arises naturally when studying Newton-like
problems of the body of least resistance. The first problem of such
kind was considered by Newton itself \cite{N}. Recently there were
made several works concerning the problem of least resistance in
various classes of admissible bodies; see, e.\,g.,
\cite{BK}-\cite{P2}, \cite{temperature}. The solution of a
minimization problem for the case of {\it rotating bodies} can be
naturally identified with a rough body (\cite{sb-math-averaged04};
see also concluding remarks to this paper).

There are many papers on particle scattering by rough bodies (see,
e.\,g., \cite{Bar}-\cite{AK-rough}); they describe bodies and flows
of particles that occur in nature. On the contrary, we assume that a
rough structure can be "manufactured", and our aim is to describe
all possible rough structures.

\subsection{Definition of a rough body}

It is supposed that the "microscopic structure"{} of the boundary of
a rough body can be detected from observations of particle
scattering on the body. From this point of view, two rough bodies
are considered equal if they scatter flows of particles in an
identical manner. Having these observations in mind, we give the
definition of a rough body.

Let $B$ be a bounded convex body. Denote by $n(\xi)$ the unit outer
normal vector to $\pl B$ at a regular point $\xi \in \pl B$, and
denote by $\set$ the set of pairs $(\xi, v) \in \pl B \times S^1$
such that $\langle n(\xi),\, v \rangle \ge 0$. Here and in what
follows, $\langle \, \cdot, \cdot \rangle$ means the standard scalar
product in $\RRR^2$. The set $\set$ is equipped with the measure
$\mu$ which is defined by $d\mu(\xi, v) = \langle n(\xi),\, v
\rangle\, d\xi\, dv$, where $d\xi$ and $dv$ are the one-dimensional
Lebesgue measures on $\pl B$ and $S^1$, respectively.

Let $\QQQ$ be a set with piecewise smooth boundary contained in $B$;
consider the billiard in $\RRR^2 \setminus \QQQ$. Note that $\QQQ$
is not necessarily connected. For $(\xi, v) \in \set$, consider a
billiard particle starting at the point $\xi$ with the velocity
$-v$. After several (maybe none) reflections from $\pl\QQQ \setminus
\pl B$, the particle will intersect $\pl B$ again, at a point $\xi^+
= \xi^+_{\QQQ,B} (\xi, v) \in \pl B$; denote by $v^+ = v^+_{\QQQ,B}
(\xi, v)$ the velocity at this point. It may happen that the initial
point $\xi$ belongs to $\pl\QQQ$; in that case we have $\xi^+ = \xi$
and the vector $v^+$ is symmetric to $v$ with respect to $n(\xi)$.
It may also happen that at some moment the particle either gets into
a singular point of $\pl\QQQ$, or touches $\pl\QQQ$ at a regular
point, or stays in $B \setminus \QQQ$ forever and does not intersect
$\pl B$ again, or makes an infinite number of reflections in finite
time. The set of corresponding points $(\xi, v)$ has zero measure,
and the corresponding values $\xi^+_{\QQQ,B} (\xi, v)$
and$v^+_{\QQQ,B} (\xi, v)$ are not defined.
 
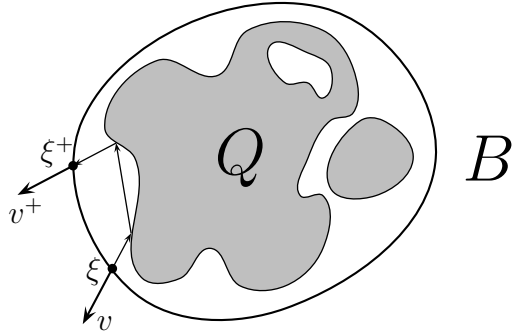
\vspace{0mm}

 \vspace*{45mm}

 \rput(5,0){
 \psecurve(-0.7,2.5)(0,0.2)(2,0)(4,2)(2.5,4)(-0.7,2.5)(0,0.2)(2,0)
 \psecurve[linewidth=0.5pt,fillcolor=lightgray,fillstyle=solid]
  (3.0,2.5)(2.6,1.6)(3.2,1.5)(3.7,2)(3.3,2.5)(3.0,2.5)(2.6,1.6)(3.2,1.5)
  \psecurve[linewidth=0.5pt,fillcolor=lightgray,fillstyle=solid](-0.1,2)(0.1,0.4)
  (0.6,0.4)(1,0.8)(1.4,0.3)(2,0.3)(2.5,0.8)(2.5,1.5)
 (2.2,1.5)(2.3,2)(2.6,2.5)(2.9,2.6)(2.9,3.2)(2.5,3.7)(1.5,3.7)(1.1,3.1)(0.3,3.3)
 (-0.4,2.4)(0,2)(0.1,0.4)(0.6,0.4)
  \psecurve[linewidth=0.5pt,fillcolor=white,fillstyle=solid](2.0,3.0)(2.3,3.0)
  (2.5,2.8)(2.6,3.2)(2.3,3.4)(1.8,3.5)(2.0,3.0)(2.3,3.0)(2.5,2.8)
 \rput(1.4,2){\Huge $\QQQ$}
 \rput(4.7,2){\Huge $B$}
 \psdots(-0.3,0.53)(-0.82,1.9)
 \psline[linewidth=0.5pt]{->}(-0.3,0.53)(-0.05,1)
  \psline[linewidth=0.5pt]{->}(-0.05,1)(-0.25,2.2)
    \psline[linewidth=0.5pt]{->}(-0.25,2.2)(-0.82,1.9)
 \psline[linewidth=0.8pt,arrows=->,arrowscale=1.5](-0.3,0.53)(-0.69325,-0.20931)
 \psline[linewidth=0.8pt,arrows=->,arrowscale=1.5](-0.82,1.9)(-1.561,1.51)
 \rput(-0.4,-0.2){$v$}
 \rput(-0.55,0.5){$\xi$}
 \rput(-1.45,1.3){$v^+$}
 \rput(-1.03,2.08){$\xi^+$}
 }
  \vspace{10mm}

Thus, there is defined the one-to-one mapping $T_{\QQQ,B}:\, (\xi,
v) \mapsto (\xi^+_{\QQQ,B}(\xi, v),\, v^+_{\QQQ,B}(\xi, v))$ of a
full measure subset of $\set$ onto itself. It has the following
properties:
 \begin{quote}
{\bf T1.}~ $T_{\QQQ,B}$ preserves the measure $\mu$.

{\bf T2.}~ $T_{\QQQ,B}^{-1} = T_{\QQQ,B}$.
 \end{quote}

The mapping $T_{\QQQ,B}$ induces the measure $\nu_{\QQQ,B}$ on
$\mathbb T^3 = S^1 \times S^1 \times S^1$ in the following way. Let
$A \subset \mathbb T^3$ be a Borel set; by definition,
$$
\nu_{\QQQ,B}(A) = \mu \left( \{ (\xi, v) \in \set:\, (v,\,
v^+_{\QQQ,B} (\xi, v), n(\xi)) \in A \} \right).
$$
In fact, the measure $\nu_{\QQQ,B}$ contains information about
particle scattering on $\QQQ$. Imagine that an observer has no means
to track the trajectory of particles inside $B$. Instead, for each
incident particle there is registered the triple of vectors: the
initial and final velocities (measured at the points of first and
second intersection with $\pl B$), and the normal vector to $\pl B$
at the point of first intersection with $\pl B$. The normal vector
at the second point of intersection is not registered; as will be
seen later on (lemma 1), if the area of $B \setminus \QQQ$ is small
then the difference between the normal vectors at these two points
is also small. The measure $\nu_{Q,B}$ describes the distribution of
triples.

\begin{opr}\label{o1}
We say that a sequence of sets $\{ \QQQ_m,\, m = 1,\, 2,\ldots \}$
represents a {\rm rough body}, if
 \begin{quote}
{\bf M1.}~ $\QQQ_m \subset B$ and Area$(B \setminus \QQQ_m) \to 0$
as $m \to \infty$;

{\bf M2.}~ the sequence of measures $\nu_{\QQQ_m,B}$ weakly
converges.
 \end{quote}
Two sequences of such sets are called equivalent, if the
corresponding limiting measures coincide. An equivalence class is
called a {\rm body obtained by roughening} $B$, or simply {\rm rough
body}, and denoted by $\BBB$, and the corresponding limiting measure
is denoted by $\nu_\BBB$.
\end{opr}

Note that the sets $Q_m$ in this definition are not necessarily
connected.
 \vspace{2mm}

{\bf Remark.} {\it Since $\TTT^3$ is compact and the full measure of
$\TTT^3$ satisfies $\nu_{\QQQ,B}(\TTT^3) \le 2\pi\, |\pl B|$, one
concludes that the set of measures $\{ \nu_{\QQQ,B} \}$, with fixed
$B$, is weakly precompact. That is, any sequence of measures $\{
\nu_{\QQQ_m,B} \}$ contains a weakly converging subsequence. In this
sense one can say that a sequence, satisfying only the condition M1,
can represent more than one rough body.}
 \vspace{1mm}

We would also like to mention that, firstly, two rough bodies
obtained one from another by translation are identified, according
to our definition. Secondly, particle scattering on $\BBB$ in a
small neighborhood of $\xi \in \pl B$ can be detected if $\xi$ is an
extreme point of $B$, and cannot otherwise. Indeed, if $\xi$ is an
extreme point of $B$, the scattering is described by the restriction
of $\nu_\BBB$ on $\TTT^2 \times \NNN_{n(\xi)}$, with $\NNN_{n(\xi)}$
being a small neighborhood of $n(\xi)$ in $S^1$. If, otherwise,
$\xi$ is not an extreme point of $B$, that is, belongs to an open
linear segment contained in $\pl B$, the scattering can only be
determined on the whole segment.

Actually, from the viewpoint of applications to the problems of
optimal resistance in {\it homogeneous} and {\it rarefied} media
(see the section Concluding Remarks and Applications), these
drawbacks are not so serious. Indeed, resistance of a body is
invariant under translations (due to homogeneity). Besides, if the
boundary of a body contains a linear segment, one does not need to
know scattering at each point of the segment; it suffices to know it
on the segment in the whole (due to homogeneity and rarefaction).

The definition of a rough body could be made in a slightly different
way, basing on measures defined on $S^1 \times S^1 \times \pl B$. In
that case the triple $(v, v^+, \xi)$ should be registered, with
$\xi$ being the point of first intersection with $\pl B$. That
definition would allow one to register particle scattering at each
point of $\pl B$ and to distinguish between bodies obtained by
translation one from another. However, we prefer to adopt the former
definition, since it seems to us mathematically more transparent and
makes the arguments a bit easier.

\subsection{Examples}

Sometimes it is convenient to use another representation of the
measure $\nu_\BBB$. Namely, consider the change of coordinates $(v,
v^+, n) \mapsto (\vphi, \vphi^+, n)$, where $\vphi = \Arg\, v -
\Arg\, n$,\, $\vphi^+ = \Arg\, v^+ - \Arg\, n$. Here $\Arg\, v$ is
the angle between a fixed vector and $v$ measured, say, clockwise
from this vector to $v$. If $(v, v^+, n) \in$ spt\,$\nu_\BBB$ then
$\vphi$ and $\vphi^+$ belong to $\interval$ modulo $2\pi$. Introduce
the shorthand notation $\Box := \interval \times \interval$ and
define the mapping $\varpi: \Box \times S^1 \to \TTT^3$ by
$\varpi(\vphi, \vphi^+, n) = (v, v^+, n)$. One has spt\,$\nu_\BBB
\subset \varpi(\Box \times S^1)$. Denote $\breve\nu_\BBB :=
(\varpi^{-1})^\# \nu_\BBB$. Sometimes this measure can be
factorized: $\breve\nu_\BBB = \eta_\BBB \otimes \tau_B$, where
$\eta_\BBB$ is defined on $\Box$ and $\tau_B$ is the surface measure
on $B$; so to say, the "roughness"{} is "homogeneous"{} along the
body's boundary. Consider several examples.
 \vspace{1mm}

{\bf Example 1}:~ {\it "smooth body".}\\
The rough body represented by the sequence $Q_m = B$ is identified
with $B$ itself. The corresponding measure is $\breve\nu_B = \eta_0
\otimes \tau_B$, where the measure $\eta_0$ has the density
$\cos\vphi \cdot \del(\vphi + \vphi^+)$; the support of $\eta_0$ is
shown on the figure below. On this figure, $B$ is taken to be an
ellipse.
  
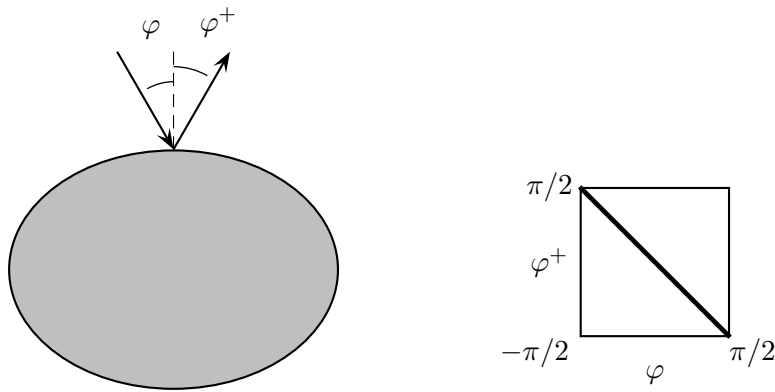
\vspace{46mm}

 \rput(-0.4,0){
 \psellipse[fillstyle=solid,fillcolor=lightgray](3,1.2)(2.2,1.6)
 \psline[arrows=->,arrowscale=1.8](2.25,4.1)(3,2.8)
 \psline[arrows=<-,arrowscale=1.8](3.75,4.1)(3,2.8)
 \psline[linewidth=0.4pt,linestyle=dashed](3,4.1)(3,2.8)
  \psarc[linewidth=0.5pt](3,3){0.9}{60}{90}
  \psarc[linewidth=0.5pt](3,3){0.7}{90}{120}
 \rput(2.7,4.4){\small $\vphi$}
 \rput(3.6,4.5){\small $\vphi^+$}
  }
 \rput(8,0.3){
 \psframe(0,0)(2,2)
 \psline[linewidth=1.8pt](0,2)(2,0)
 \rput(1,-0.5){\small $\vphi$}
 \rput(-0.4,1){\small $\vphi^+$}
  \rput(-0.6,-0.2){\small $-\pi/2$}
  \rput(2.3,-0.2){\small $\pi/2$}
  \rput(-0.4,2){\small $\pi/2$}
  }

\vspace{11mm}

{\bf Example 2}:~ {\it roughness formed by triangular hollows}.\\
$Q_m$ is a $2m$-polygon; the angles $270^0$ alternate with the
angles that are slightly smaller than $90^0$. All vertices
corresponding to the angles smaller than $90^0$ belong to $\pl B$.
Any two sides that form an angle $270^0$ are equal. The largest side
length tends to zero as $m \to \infty$. Thus, the set $Q_m$ is
obtained by moving off $m$ "hollows"{} from its convex hull, each of
the hollows being an isosceles right triangle.

The corresponding measure is $\breve\nu_\BBB = \eta_\triangledown
\otimes \tau_B$, where the measure $\eta_\triangledown$ has the
density $$\cos\vphi \cdot \left[ \chi_{[-\pi/2,-\pi/4]}(\vphi)\,
\del(\vphi + \vphi^+ + \frac \pi 2) + \chi_{[-\pi/4,\pi/4]}(\vphi)\,
\del(\vphi - \vphi^+) + \right.
 $$
 $$
\left. + \chi_{[\pi/4,\pi/2]}(\vphi)\, \del(\vphi + \vphi^+ - \frac
\pi 2) \right] + |\sin\vphi| \cdot \left[ \chi_{[-\pi/4,0]}(\vphi)\,
\del(\vphi + \vphi^+ + \frac \pi 2) - \right.
 $$
 $$
\left. - \chi_{[-\pi/4,\pi/4]}(\vphi)\, \del(\vphi - \vphi^+) +
\chi_{[0,\pi/4]}(\vphi)\, \del(\vphi + \vphi^+ - \frac \pi 2)
\right].$$
 Thus, the support of $\eta_\triangledown$ is the union of
three segments; see the figure below. The middle segment $\vphi^+ =
\vphi$ corresponds to double reflections, and the lateral segments,
$\vphi^+ = -\vphi - \pi/2$ and $\vphi^+ = -\vphi + \pi/2$,
correspond to single reflections, from the right or from the left
side of a triangular hollow.  On the figure, $B$ is a circle.
 
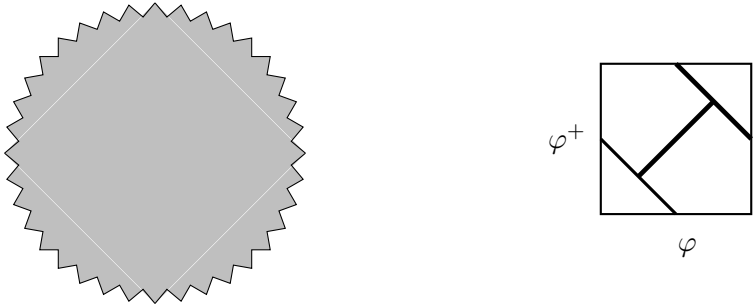
\vspace{0mm}

  \vspace*{48mm}

 \rput(3,2.5){
 \scalebox{2}{
                \pspolygon[linewidth=0pt,linecolor=white,fillstyle=solid,fillcolor=lightgray](1,0)(0,1)(-1,0)(0,-1)
      \psline[linewidth=0.2pt,fillstyle=solid,fillcolor=lightgray](1,0)(0.9056,0.0792)(0.9848,0.17365)(0.878,0.23525)
      (0.9397,0.342)(0.82386,0.3842)(0.866,0.5)(0.7446,0.5214)(0.766,0.6428)(0.6428,0.6428)
      (0.6428,0.766)(0.5214,0.7446)(0.5,0.866)(0.3842,0.82386)(0.342,0.9397)(0.23525,0.878)
      (0.17365,0.9848)(0.0792,0.9056)(0,1)
               \psline[linewidth=0.2pt,fillstyle=solid,fillcolor=lightgray](-1,0)(-0.9056,0.0792)(-0.9848,0.17365)(-0.878,0.23525)
      (-0.9397,0.342)(-0.82386,0.3842)(-0.866,0.5)(-0.7446,0.5214)(-0.766,0.6428)(-0.6428,0.6428)
      (-0.6428,0.766)(-0.5214,0.7446)(-0.5,0.866)(-0.3842,0.82386)(-0.342,0.9397)(-0.23525,0.878)
      (-0.17365,0.9848)(-0.0792,0.9056)(0,1)
            \psline[linewidth=0.2pt,fillstyle=solid,fillcolor=lightgray](1,0)(0.9056,-0.0792)(0.9848,-0.17365)(0.878,-0.23525)
      (0.9397,-0.342)(0.82386,-0.3842)(0.866,-0.5)(0.7446,-0.5214)(0.766,-0.6428)(0.6428,-0.6428)
      (0.6428,-0.766)(0.5214,-0.7446)(0.5,-0.866)(0.3842,-0.82386)(0.342,-0.9397)(0.23525,-0.878)
      (0.17365,-0.9848)(0.0792,-0.9056)(0,-1)
               \psline[linewidth=0.2pt,fillstyle=solid,fillcolor=lightgray](-1,0)(-0.9056,-0.0792)(-0.9848,-0.17365)(-0.878,-0.23525)
      (-0.9397,-0.342)(-0.82386,-0.3842)(-0.866,-0.5)(-0.7446,-0.5214)(-0.766,-0.6428)(-0.6428,-0.6428)
      (-0.6428,-0.766)(-0.5214,-0.7446)(-0.5,-0.866)(-0.3842,-0.82386)(-0.342,-0.9397)(-0.23525,-0.878)
      (-0.17365,-0.9848)(-0.0792,-0.9056)(0,-1)
       }
       }

        \rput(9,4.2){
 \pspolygon[linewidth=0.8pt](0,0)(2,0)(2,-2)(0,-2)
  \psline[linewidth=1.2pt](0,-1)(1,-2)
   \psline[linewidth=1.8pt](2,-1)(1,0)
     \psline[linewidth=1.8pt](0.5,-1.5)(1.5,-0.5)
       \rput(-0.45,-1){\small $\vphi^+$}
    \rput(1.15,-2.45){\small $\vphi$}
     }

  \vspace*{0mm}

{\bf Example 3}:~ {\it roughness formed by rectangular hollows}.\\
The sets $Q_m$ are obtained by removing a finite number of
"rectangular hollows"{} from $B$. In other words, one has $Q_m = B
\setminus \left( \cup_n \Om_{m,n} \right)$, where the removed sets
$\Om_{m,n}$ do not mutually intersect and each set $\pl\Om_{m,n}
\setminus \pl B$ is the union of three sides of a rectangle. The
ratio (width)/(depth) of a hollow depends only on $m$ and is denoted
by $h_m$. Denote by $l_m = |\pl B \setminus \cup_n \left(
\pl\Om_{m,n} \right)| / |\pl B|$ the relative length of the part of
boundary $\pl B$ not covered by hollows. We assume that
$\lim_{m\to\infty} h_m = 0 = \lim_{m\to\infty} l_m$. On the figure
below, $B$ is a square.
 
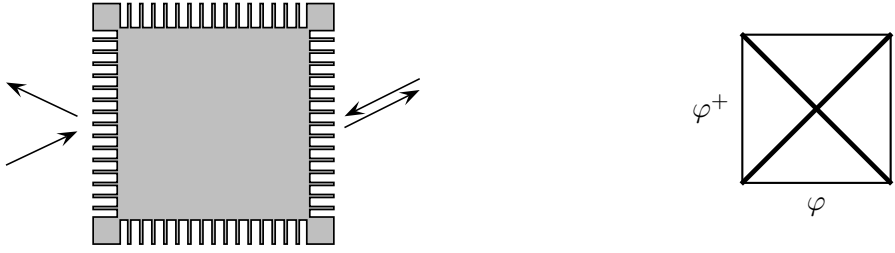
\vspace{0mm}

 \vspace*{40mm}

 \rput(3,2){
       \psline[linewidth=0.6pt,arrows=->,arrowscale=2](-2.65,-0.55)(-1.7,-0.1)
        \psline[linewidth=0.6pt,arrows=<-,arrowscale=2](-2.65,0.55)(-1.7,0.1)
       \psline[linewidth=0.6pt,arrows=->,arrowscale=2](2.85,0.6)(1.85,0.1)
        \psline[linewidth=0.6pt,arrows=<-,arrowscale=2](2.85,0.45)(1.85,-0.05)
 \scalebox{0.8}{
 \pspolygon[fillstyle=solid,fillcolor=lightgray]     (-1.55,-1.6)(-1.55,-2)(-2,-2)(-2,-1.55)(-1.6,-1.55)
                                                     (-1.6,-1.425)(-2,-1.425)(-2,-1.375)(-1.6,-1.375)
                                                     (-1.6,-1.225)(-2,-1.225)(-2,-1.175)(-1.6,-1.175)
                                                     (-1.6,-1.025)(-2,-1.025)(-2,-0.975)(-1.6,-0.975)
                                                     (-1.6,-0.825)(-2,-0.825)(-2,-0.775)(-1.6,-0.775)
                                                     (-1.6,-0.625)(-2,-0.625)(-2,-0.575)(-1.6,-0.575)
                                                     (-1.6,-0.425)(-2,-0.425)(-2,-0.375)(-1.6,-0.375)
                                                     (-1.6,-0.225)(-2,-0.225)(-2,-0.175)(-1.6,-0.175)
                                                     (-1.6,-0.025)(-2,-0.025)(-2,0.025)(-1.6,0.025)
                                                     (-1.6,0.175)(-2,0.175)(-2,0.225)(-1.6,0.225)
                                                     (-1.6,0.375)(-2,0.375)(-2,0.425)(-1.6,0.425)
                                                     (-1.6,0.575)(-2,0.575)(-2,0.625)(-1.6,0.625)
                                                     (-1.6,0.775)(-2,0.775)(-2,0.825)(-1.6,0.825)
                                                     (-1.6,0.975)(-2,0.975)(-2,1.025)(-1.6,1.025)
                                                     (-1.6,1.175)(-2,1.175)(-2,1.225)(-1.6,1.225)
                                                     (-1.6,1.375)(-2,1.375)(-2,1.425)(-1.6,1.425)
                                                     (-1.6,1.55)(-2,1.55)(-2,2)(-1.55,2)(-1.55,1.6)
                                                     (-1.425,1.6)(-1.425,2)(-1.375,2)(-1.375,1.6)
                                                     (-1.225,1.6)(-1.225,2)(-1.175,2)(-1.175,1.6)
                                                     (-1.025,1.6)(-1.025,2)(-0.975,2)(-0.975,1.6)
                                                     (-0.825,1.6)(-0.825,2)(-0.775,2)(-0.775,1.6)
                                                     (-0.625,1.6)(-0.625,2)(-0.575,2)(-0.575,1.6)
                                                     (-0.425,1.6)(-0.425,2)(-0.375,2)(-0.375,1.6)
                                                     (-0.225,1.6)(-0.225,2)(-0.175,2)(-0.175,1.6)
                                                     (-0.025,1.6)(-0.025,2)(0.025,2)(0.025,1.6)
                                                     (0.175,1.6)(0.175,2)(0.225,2)(0.225,1.6)
                                                     (0.375,1.6)(0.375,2)(0.425,2)(0.425,1.6)
                                                     (0.575,1.6)(0.575,2)(0.625,2)(0.625,1.6)
                                                     (0.775,1.6)(0.775,2)(0.825,2)(0.825,1.6)
                                                     (0.975,1.6)(0.975,2)(1.025,2)(1.025,1.6)
                                                     (1.175,1.6)(1.175,2)(1.225,2)(1.225,1.6)
                                                     (1.375,1.6)(1.375,2)(1.425,2)(1.425,1.6)
                                                     (1.55,1.6)(1.55,2)(2,2)
                                                                           (2,1.55)(1.6,1.55)
                                                     (1.6,1.425)(2,1.425)(2,1.375)(1.6,1.375)
                                                     (1.6,1.225)(2,1.225)(2,1.175)(1.6,1.175)
                                                     (1.6,1.025)(2,1.025)(2,0.975)(1.6,0.975)
                                                     (1.6,0.825)(2,0.825)(2,0.775)(1.6,0.775)
                                                     (1.6,0.625)(2,0.625)(2,0.575)(1.6,0.575)
                                                     (1.6,0.425)(2,0.425)(2,0.375)(1.6,0.375)
                                                     (1.6,0.225)(2,0.225)(2,0.175)(1.6,0.175)
                                                     (1.6,0.025)(2,0.025)(2,-0.025)(1.6,-0.025)
                                                     (1.6,-0.175)(2,-0.175)(2,-0.225)(1.6,-0.225)
                                                     (1.6,-0.375)(2,-0.375)(2,-0.425)(1.6,-0.425)
                                                     (1.6,-0.575)(2,-0.575)(2,-0.625)(1.6,-0.625)
                                                     (1.6,-0.775)(2,-0.775)(2,-0.825)(1.6,-0.825)
                                                     (1.6,-0.975)(2,-0.975)(2,-1.025)(1.6,-1.025)
                                                     (1.6,-1.175)(2,-1.175)(2,-1.225)(1.6,-1.225)
                                                     (1.6,-1.375)(2,-1.375)(2,-1.425)(1.6,-1.425)
                                                     (1.6,-1.55)(2,-1.55)(2,-2)(1.55,-2)(1.55,-1.6)
                                                     (1.425,-1.6)(1.425,-2)(1.375,-2)(1.375,-1.6)
                                                     (1.225,-1.6)(1.225,-2)(1.175,-2)(1.175,-1.6)
                                                     (1.025,-1.6)(1.025,-2)(0.975,-2)(0.975,-1.6)
                                                     (0.825,-1.6)(0.825,-2)(0.775,-2)(0.775,-1.6)
                                                     (0.625,-1.6)(0.625,-2)(0.575,-2)(0.575,-1.6)
                                                     (0.425,-1.6)(0.425,-2)(0.375,-2)(0.375,-1.6)
                                                     (0.225,-1.6)(0.225,-2)(0.175,-2)(0.175,-1.6)
                                                     (0.025,-1.6)(0.025,-2)(-0.025,-2)(-0.025,-1.6)
                                                     (-0.175,-1.6)(-0.175,-2)(-0.225,-2)(-0.225,-1.6)
                                                     (-0.375,-1.6)(-0.375,-2)(-0.425,-2)(-0.425,-1.6)
                                                     (-0.575,-1.6)(-0.575,-2)(-0.625,-2)(-0.625,-1.6)
                                                     (-0.775,-1.6)(-0.775,-2)(-0.825,-2)(-0.825,-1.6)
                                                     (-0.975,-1.6)(-0.975,-2)(-1.025,-2)(-1.025,-1.6)
                                                     (-1.175,-1.6)(-1.175,-2)(-1.225,-2)(-1.225,-1.6)
                                                     (-1.375,-1.6)(-1.375,-2)(-1.425,-2)(-1.425,-1.6)

                                                      } }
 \rput(10,1.2){
 \psframe(0,0)(2,2)
 \psline[linewidth=1.8pt](0,2)(2,0)
  \psline[linewidth=1.8pt](0,0)(2,2)
 \rput(1,-0.3){\small $\vphi$}
 \rput(-0.4,1){\small $\vphi^+$}
  }

  \vspace{3mm}

The measure $\breve\nu_\BBB$ equals $\breve\nu_\BBB = \eta_\square
\otimes \tau_B$. The density of the measure $\eta_\square$ equals
$\frac 12\, \cos\vphi \cdot (\del(\vphi + \vphi^+) + \del(\vphi -
\vphi^+))$, and the support is the union of two diagonals, $\vphi^+
= \vphi$ and $\vphi^+ = -\vphi$; see the figure. The particles with
even (odd) number of reflections contribute to the first (second)
diagonal.

\subsection{Main theorem}

According to the definition \ref{o1}, each rough body is identified
with a measure on $\mathbb T^3$. The question is: what is the set of
these measures? The following definition and theorem give the
answer.

Let us first introduce some notation: $\pi_{v,n}: \TTT^3 \to
\TTT^2$,\, $\pi_n: \TTT^3 \to S^1$, etc. are projections onto the
corresponding subspaces: $\pi_{v, n}(v, v^+, n) = (v, n)$,\,
$\pi_{n}(v, v^+, n) = n$, etc.; $\pi_d: \TTT^3 \to \TTT^3$ is the
symmetry with respect to the plane $v = v^+$, that is, $\pi_d(v,
v^+, n) = (v^+, v, n)$;\, $z_+ = \max \{ 0, z \}$ is the positive
part of $z \in \RRR$; and $u$ means Lebesgue measure on $S^1$.
Recall that $\tau_B$ is the surface measure on $B$ and is defined on
$S^1$.

\begin{opr}\label{oM}
We denote by $\MMM_B$ the set of measures $\nu$ on $\mathbb T^3$
such that
\begin{quote}
{\bf A1}~ the marginal measures $\pi_{v,n}^\# \nu$ and
$\pi_{v^+,n}^\# \nu$ are
$$
\pi_{v,n}^\# \nu\ =\ \langle v, n \rangle_+ \cdot u \otimes
\tau_B,~~~~~ \pi_{v^+,n}^\# \nu\ =\ \langle v^+, n \rangle_+ \cdot u
\otimes \tau_B;~~~~~~~~~~
$$
$$
\text{\bf A2} \hspace{5mm} \pi_d^\# \nu = \nu. \hspace{120mm}
$$
 \end{quote}
Denote also $\MMM = \cup_B \MMM_B$, the union being taken over all
bounded convex bodies $B$.
\end{opr}

Taking into account the Alexandrov theorem on characterization of
surface measures, one concludes that $\MMM$ is the set of measures
$\nu$ on $\TTT^3$ such that

1) the marginal measure $\pi_n^\# \nu =: \tau$ satisfies the
conditions

\hspace*{5mm} 1a.~ $\int_{S^1} n\, d\tau(n) = 0$;

\hspace*{5mm} 1b.~ for any $v \in S^1$ holds $\int_{S^1} \langle n,
v \rangle^2 \, d\tau(n) \ne 0$;

2) the marginal measures $\pi_{v,n}^\# \nu$ and $\pi_{v^+,n}^\# \nu$
satisfy the conditions

\hspace*{5mm} 2a.~ $\pi_{v,n}^\# \nu = \langle v, n \rangle_+ \cdot
u \otimes \tau$;

\hspace*{5mm} 2b.~ $\pi_{v^+,n}^\# \nu = \langle v^+, n \rangle_+
\cdot u \otimes \tau$.

Thus, these marginal measures coincide; the only difference is in
the notation for the variables: $v,n$ in the case 2a and $v^+,n$ in
the case 2b.
 \vspace{2mm}

Now we can state the main theorem.
 \vspace{2mm}

{\bf {Theorem}.} {\it The set of measures $\{ \nu_\BBB \}$, with
$\BBB$ being all possible bodies obtained by roughening $B$,
coincides with $\MMM_B$. Therefore, $\left\{ \nu_\BBB,\ \BBB\,
\right.$is a rough body\,$\left.\! \right\} = \MMM$.}
 \vspace{2mm}

In section 2, we formulate two auxiliary lemmas and using them,
prove the theorem. In section 3, the lemmas are proved. Section 4
contains concluding remarks and applications of theorem to problems
of optimal aerodynamic resistance. Appendices A and B contain proofs
of some auxiliary technical results.

\section{Statement of auxiliary lemmas\\ and proof of theorem}

\subsection{Statement of lemma 1}

Fix a bounded convex body $B$. Two points $\xi_1$,\, $\xi_2 \in \pl
B$,\, $\xi_1 \ne \xi_2$ divide the curve $\pl B$ into two arcs.
Denote by $l(\xi_1, \xi_2)$ the length of the smallest arc and
denote
$$
c = c_B := \inf_{ \scriptsize \begin{array}{c} \xi_1, \xi_2 \in \pl
B\\ \xi_1 \ne \xi_2
\end{array}} \frac{|\xi_1 - \xi_2|}{l(\xi_1, \xi_2)}\,;
$$
one obviously has $0 < c < 1$.

Let $\QQQ \subset B$;  denote
\begin{equation*}\label{xi-xi+}
\overline{|\xi - \xi^+|}_{\QQQ,B} := \int\!\!\!\int_{\set} |\xi -
\xi^+_{\QQQ,B}(\xi, v)|\, d\mu(\xi,v)
\end{equation*}
and
\begin{equation*}\label{Del theta}
\overline{|n - n^+|}_{\QQQ,B} := \int\!\!\!\int_{\set} |n(\xi) -
n(\xi^+_{\QQQ,B}(\xi, v))|\, d\mu(\xi,v).
\end{equation*}

\begin{lem}\label{l1}
(a) The following holds true:
$$
\overline{|\xi - \xi^+|}_{\QQQ,B} \le 2\pi \cdot \text{\rm Area} (B
\setminus \QQQ).
$$

\hspace*{13mm} (b) For sufficiently small {\rm Area}$(B \setminus
\QQQ)$,\footnote{That is, smaller that a value depending only on
$B$.} one has
$$
\overline{|n - n^+|}_{\QQQ,B} \le \frac{2\pi
\sqrt{8\pi}}{\sqrt{c}}\, \sqrt{\text{\rm Area} (B \setminus \QQQ)}.
$$
\end{lem}

\subsection{Statement of lemma 2}

Let us first introduce the notion of a hollow.

\begin{opr}\label{o3}
Let $\Om \subset \RRR^2$ be a closed bounded set with piecewise
smooth boundary and $I \subset \pl \Om$, where

(i)~ $I$ is an interval contained in a straight line $\langle x, n
\rangle = a$ and

(ii)~ $\Om \setminus I$ is contained in the open half-plane $\langle
x, n \rangle < a$. Here $n$ is a fixed unit vector.\\
Then the pair $(\Om, I)$ is called a hollow oriented by $n$, or just
an $n$-hollow.
\end{opr}

Here and in what follows, $I$ is shown dashed, and $\pl \Om
\setminus I$ is shown by solid line.
 \vspace{0mm}

 \newpage

 
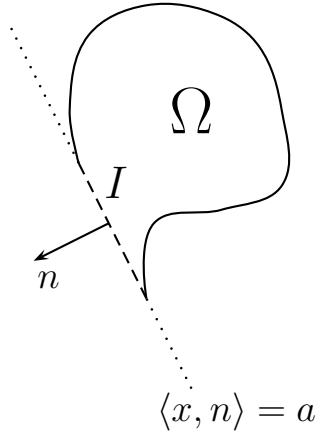
\vspace*{35mm}

 \rput(5.5,1.4){
 \psline[linestyle=dashed](-0.4,0.8)(0.5,-1)
 \psline[linestyle=dotted,linewidth=1pt](-1.3,2.6)(-0.4,0.8)
 \psline[linestyle=dotted,linewidth=1pt](0.5,-1)(1.1,-2.2)
 \psline[arrows=->,arrowscale=1.3](0,0)(-1,-0.5)
 \pscurve(-0.4,0.8)(-0.52,1.5)(-0.27,2.5)(0.5,2.9)(1.5,2.75)(2,2.35)(2.3,1.5)(2.3,0.5)
 (1.5,0.18)(0.6,0)(0.5,-1)
 \rput(0.1,0.55){\Large $I$}
 \rput(1.1,1.5){\Huge $\Om$}
 \rput(-0.8,-0.75){\large $n$}
 \rput(1.7,-2.5){\large $\langle x,n \rangle = a$}

 }

\vspace{22mm}

Define the measure $\tilde\mu_I$ on $I \times S^1$ by
$d\tilde\mu_I(\xi,v) = \frac{\langle n, v \rangle_+}{|I|}\, d\xi\,
dv$, where $|I|$ means the length of $I$. Obviously, $\tilde\mu_I$
is supported on the set $\iset := \{ (\xi, v) \in I \times S^1 :
\langle n, v \rangle \ge 0 \}$. Define the one-to-one mapping $(\xi,
v) \mapsto (\Xi^+_{\Om,I}(\xi,v),\, V^+_{\Om,I}(\xi,v))$ of a full
measure subset of $\iset$ onto itself. Namely, consider the billiard
in $\Om$. Let $(\xi, v) \in \iset$; consider the billiard particle
starting at the point $\xi$ with the velocity $-v$. It makes several
reflections from $\pl\Om \setminus I$ and then reflects from $I$
again, at a point $\Xi^+ = \Xi^+_{\Om,I}(\xi,v)$. The velocity
immediately before this reflection is denoted by $V^+ =
V^+_{\Om,I}(\xi,v)$. The mapping so defined preserves the measure
$\tilde\mu_I$ and is an involution, that is, coincides with its
inverse.

One can give an equivalent definition based on the mapping
$\xi^+_{\QQQ,B}(\xi,v),\, v^+_{\QQQ,B}(\xi,v)$ just defined in
subsection 1.2. Take a set $\QQQ$ such that $\Om$ is a connected
component of conv$\QQQ \setminus \QQQ$ and $I$ is a connected
component of $\pl(\text{conv}\, \QQQ) \setminus \pl\QQQ$. For $(\xi,
v) \in \iset$, let by definition $(\Xi^+_{\Om,I}(\xi,v),\,
V^+_{\Om,I}(\xi,v)) := (\xi^+_{\QQQ, \text{conv}\QQQ}(\xi,v),\,
v^+_{\QQQ, \text{conv}\QQQ}(\xi,v))$. This definition does not
depend on the choice of $Q$.

\begin{opr}\label{o4}
Let $(\Om, I)$ be a hollow. The measure $\eta_{\Om,I}$ on $\mathbb
T^2 = S^1 \times S^1$ is defined as follows. For a Borel set $A
\subset \TTT^2$, put
$$
\eta_{\Om,I}(A) := \tilde\mu_I (\{ (\xi,v) \in \iset:\ (v,
V^+_{\Om,I}(\xi,v)) \in A \}).
$$
We shall say that $\eta_{\Om,I}$ is the measure generated by the
hollow $(\Om, I)$.
\end{opr}

Here we use the notation $\pi_{v},\, \pi_{v^+} : \TTT^2 \to S^1$ for
the projections onto the subspaces $\{ v \}$ and $\{ v^+ \}$,
respectively; $\pi_{v}(v, v^+) = v$,\, $\pi_{v^+}(v, v^+) = v^+$. We
also denote by $\pi_d$ the symmetry with respect to the diagonal $v
= v^+$;\, $\pi_d(v, v^+) = (v^+, v)$.

\begin{opr}\label{oLam} Denote by $\Lam_n$ the set of measures $\eta$ on $\TTT^2$ such that

1) $d\pi_v^\# \eta(v) = \langle v,n \rangle_+\, dv$,~ $\,
d\pi_{v^+}^\# \eta(v^+) = \langle v^+,n \rangle_+\, dv^+$;

2) $\pi_d^\# \eta = \eta$.
\end{opr}

Any measure $\eta_{\Om,I}$ generated by an $n$-hollow belongs to
$\Lam_n$. Indeed, for any $A \subset S^1$ one has $\pi_v^\#
\eta_{\Om,I}(A) = \eta_{\Om,I}(A \times S^1) = \tilde\mu_I(\{
(\xi,v) \in \iset:\ v \in A \}) = \frac{1}{|I|} \int\!\!\!\int_{I
\times A} \langle n, v \rangle_+\, d\xi\, dv = \int_A \langle n, v
\rangle_+\, dv$. This proves the first equality in 1).

Similarly, one has $\pi_{v^+}^\# \eta_{\Om,I}(A) = \eta_{\Om,I}(S^1
\times A) = \tilde\mu_I(\{ (\xi,v) \in \iset:\ V^+_{\Om,I}(\xi,v)
\in A \})$. Since the mapping $(\xi, v) \mapsto (\Xi^+_{\Om,I},
V^+_{\Om,I})$ preserves the measure, one gets the value
$\tilde\mu_I(\{ (\xi,v) \in \iset:\ v \in A \})$, which in turns
equals to $\int_A \langle n, v \rangle_+\, dv$. This proves the
second equality in 1). Finally, the relation 2) for $\eta_{\Om,I}$
is a simple consequence of involutive and measure preserving
properties of the mapping $(\xi, v) \mapsto (\Xi^+_{\Om,I},
V^+_{\Om,I})$.

\begin{lem}\label{l2}
The set of measures generated by $n$-hollows is weakly dense in
$\Lam_n$.
\end{lem}

\subsection{Proof of the direct statement of theorem}

Here we prove that for any body $\BBB$ obtained by roughening $B$
holds $\nu_\BBB \in \MMM_B$.

Let $\QQQ \subset B$; define the measure $\nu^{\,\prime}_{\QQQ,B}$
on $\mathbb T^3$ by
$$
\nu^{\,\prime}_{\QQQ,B}(A) := \mu \left( \{ (\xi, v) \in \set:\,
(v,\, v^+_{\QQQ,B} (\xi, v), n(\xi^+_{\QQQ,B}(\xi,v))) \in A \}
\right),
$$
where $A$ is an arbitrary Borel subset of $\mathbb T^3$. Thus, the
definition of both $\nu_{\QQQ,B}$ and $\nu^{\,\prime}_{\QQQ,B}$ is
based on observations of vector triples $(v, v^+, n)$ and $(v, v^+,
n^+)$, respectively. Here $n$ and $n^+$ are the outer normals to
$\pl B$ at the points where the particle {\it gets in} $B$ and {\it
gets out} of $B$. The measures $\nu_{\QQQ,B}$ and
$\nu^{\,\prime}_{\QQQ,B}$ have the following properties:
 \beq\lab{ast_1}
\pi_{v,n}^\# \nu_{\QQQ,B}\ =\ \langle v, n \rangle_+ \cdot u \otimes
\tau_B,
 \eeq
 \beq\lab{ast_2}
\pi_{v^+,n^+}^\# \nu^{\,\prime}_{\QQQ,B}\ =\ \langle v^+, n^+
\rangle_+ \cdot u \otimes \tau_B,
 \eeq
 \beq\lab{ast_3}
\pi_d^\# \nu_{\QQQ,B}\ =\ \nu^{\,\prime}_{\QQQ,B}.
 \eeq

Consider a sequence $\{ \QQQ_m \}$ representing $\BBB$; let us show
that $\nu_{\QQQ_m,B} - \nu^{\,\prime}_{\QQQ_m,B}$ weakly converges
to zero as $m \to \infty$. It is enough to prove that for any
continuous function $f$ on $\mathbb T^3$ holds
 \beq\label{ast444}
\int_{\TTT^3} f(v, v^+, n)\, d\nu_{\QQQ_m,B} (v, v^+, n) -
\int_{\TTT^3} f(v, v^+, n^+)\, d\nu^{\,\prime}_{\QQQ_m,B} (v, v^+,
n^+) \rightarrow_{m \to \infty} 0.
 \eeq
Taking into account the formulas for change of variables
$$
\int_{\TTT^3} f(v, v^+, n)\, d\nu_{\QQQ,B} (v, v^+, n) = \int_{\set}
f(v, v^+_{\QQQ,B}(\xi,v), n(\xi))\, d\mu(\xi,v)
$$
and
$$
\int_{\TTT^3} f(v, v^+, n^+)\, d\nu^{\,\prime}_{\QQQ,B} (v, v^+,
n^+) = \int_{\set} f(v, v^+_{\QQQ,B}(\xi,v),
n(\xi^+_{\QQQ,B}(\xi,v)))\, d\mu(\xi,v),
$$
the formula (\ref{ast444}) takes the form
 \beq\label{ast44}
\lim\limits_{m \to \infty} \int_{\set} \left[ f(v,
v^+_{\QQQ_m,B}(\xi,v), n(\xi^+_{\QQQ_m,B}(\xi,v))) - f(v,
v^+_{\QQQ_m,B}(\xi,v), n(\xi)) \right] d\mu(\xi,v) = 0.
 \eeq

According to lemma \ref{l1}, the difference
$n(\xi^+_{\QQQ_m,B}(\xi,v)) - n(\xi)$ converges to zero in mean,
hence it converges to zero in measure; therefore the difference
$$
f(v, v^+_{\QQQ_m,B}(\xi,v), n(\xi^+_{\QQQ_m,B}(\xi,v))) - f(v,
v^+_{\QQQ_m,B}(\xi,v), n(\xi))
$$
also converges to zero in measure. It follows that the formula
(\ref{ast44}) is true.

Thus, both $\nu_{\QQQ_m,B}$ and $\nu^{\,\prime}_{\QQQ_m,B}$ weakly
converge to $\nu_\BBB$. Substituting $\QQQ = \QQQ_m$ into the
formulas (\ref{ast_1}--\ref{ast_3}) and passing to limit as $m \to
\infty$, one gets
$$
\pi_{v,n}^\# \nu_\BBB\ =\ \langle v, n \rangle_+ \cdot u \otimes
\tau_B,
$$
$$
\pi_{v^+,n}^\# \nu_\BBB\ =\ \langle v^+, n \rangle_+ \cdot u \otimes
\tau_B,
$$
$$
\pi_d^\# \nu_\BBB\ =\ \nu_\BBB,
$$
that is, $\nu_\BBB \in \MMM_B$.

\subsection{Proof of the inverse statement of theorem}

Here it is proved that for any $\nu \in \MMM_B$ there exists a body
$\BBB$ obtained by roughening $B$ such that $\nu_\BBB = \nu$. The
proof is based on two statements.

\begin{predl}\label{p_1}
Let $B$ be a convex polygon. Then for any measure $\nu \in \MMM_B$
there exists a body $\BBB$ obtained by roughening $B$ such that
$\nu_\BBB = \nu$.
\end{predl}

\begin{proof}
Let us enumerate the sides of the polygon $B$ and denote by $c_i$
the length of the $i$th side, and by $n_i$, the outer unit normal to
this side. By $\del_n$, denote the probabilistic atomic measure on
$S^1$ concentrated at $n \in S^1$, that is, $\del_n(n) = 1$. The
surface measure of $B$ is $\tau_B = \sum c_i \del_{n_i}$; this
implies that any measure $\nu \in \MMM_B$ has the form $\nu = \sum
c_i \eta_i \otimes \del_{n_i}$, where $\eta_i \in \Lam_{n_i}$.

According to lemma \ref{l2}, any measure $\eta_i$ is the weak limit
as $m \to \infty$ of measures $\eta_{\Om_i^m, I_i^m}$ generated by a
sequence of $n_i$-hollows $(\Om_i^m, I_i^m)$. Now take a sequence of
sets $\QQQ_m$ such that conv$\, \QQQ_m = B$ and each connected
component of $B \setminus \QQQ_m$ is the image of a set $\Om_i^m$
under the composition of a homothety with positive ratio and a
translation, and additionally, the image of $I_i^m$ under this
transformation belongs to ($i$th side of $B) \setminus \pl Q_m$. We
also require that Area$(B \setminus \QQQ_m) \to 0$ and $|(i$th side
of $B) \setminus \pl Q_m| =: c_i^m \to c_i$ as $m \to \infty$. In
Appendix A it is shown how to construct such a sequence $\QQQ_m$.
The measure $\nu_{\QQQ_m, B} = \tilde \nu_m + \sum_i \nu_i^m$ is the
sum of the measure $\tilde \nu_m$ corresponding to reflections from
$\pl B \cap \pl\QQQ_m$ and the measures $\nu_i^m$ corresponding to
particles getting into the "hollows on the $i$th side". One has
$\tilde\nu_m = \sum_i (c_i - c_i^m) \cdot \eta_0 \otimes \del_{n_i}$
and $\nu_i^m = c_i^m \cdot \eta_{\Om_i^m, I_i^m} \otimes
\del_{n_i}$. The norm of $\tilde\nu_m$ goes to zero and $\nu_i^m$
weakly converges to $c_i\, \eta_i \otimes \del_{n_i}$ for any $i$;
it follows that $\nu_{\QQQ_m, B}$ weakly converges to $\nu$ as $m
\to \infty$. Therefore, the sequence $\QQQ_m$ represents a body
$\BBB$ obtained by roughening $B$, and $\nu_\BBB = \nu$.
\hspace{5mm}
\end{proof}

\begin{predl}\label{p_2}
For any measure $\nu \in \MMM_B$ there exist a sequence of convex
polygons $B_k \subset B$ with Area$(B \setminus B_k) \to 0$ and a
sequence of measures $\nu_k \in \MMM_{B_k}$ weakly converging to
$\nu$ as $k \to \infty$.
\end{predl}

\begin{proof}
Consider a partition of the circumference $S^1$ into a finite number
of arcs, $S^1 = \cup_i \SSS^{i}$. It induces the partition of $\pl
B$ into arcs $\pl B^i = \{ \xi \in \pl B : n(\xi) \in \SSS^i \}$.
Consider the polygon $\check{B}$ inscribed into $\pl B$ whose
vertices are separation points of this partition. Denote by $n_i$
the outer normal to the $i$th side of this polygon. Denote by
$s_{v_1,v_2}$ the operator of rotation on $S^1$ that takes $v_1$ to
$v_2$, and define the mapping $\Upsilon_i : \TTT^2 \times \SSS^{i}
\to \TTT^2$ by $\Upsilon_i (v, v^+, n) = (s_{n,n_i} v, s_{n,n_i}
v^+)$. Finally, consider the measure $\check{\nu} = \sum_i |b^i|\,
\eta_i \otimes \del_{n_i}$, where $|b^i|$ is the length of the $i$th
side of the polygon, and the measure $\eta_i$ on $\TTT^2$ is defined
by $\eta_i(A) = \frac{1}{|\pl B^i|}\, \nu(\Upsilon_i^{-1}(A))$ for
arbitrary Borel set $A \subset \TTT^2$. Here $|\pl B^i|$ is the
length of the arc $\pl B^i$. One easily verifies that $\check{\nu}$
belongs to $\MMM_{\check{B}}$.

Now take a sequence of partitions of $S^1$,\, $\{ \SSS_k^{i}
\}_i$,\, $k = 1,\, 2, \ldots$, where the maximum arc length of a
partition goes to zero as $k \to \infty$. Denote by $\{ \pl B_k^i
\}_i$,\, $k = 1,\, 2, \ldots$ the sequence of induced partitions of
$\pl B$, and take the sequence of polygons $B_k$ generated by these
partitions. One clearly has Area$(B \setminus B_k) \to 0$ and
\begin{equation}\label{max i}
\max_i \frac{|b_k^i|}{|\pl B_k^i|} \to 1 ~~~~~  \text{ as }~ k \to
\infty,
\end{equation}
where $|b_k^i|$ is the length of the $i$th side of $B_k$. In the
same way as above, one defines the mappings $\Upsilon_{ik} : \TTT^2
\times \SSS_k^{i} \to \TTT^2$ and the measures ${\nu_k} = \sum_i
|b_k^i|\, \eta_{ik} \otimes \del_{n_{ik}} \in \MMM_{B_k}$, where
$\eta_{ik}$ is given by $\eta_{ik}(A) := \frac{1}{|\pl B_k^i|}\,
\nu(\Upsilon_{ik}^{-1}(A))$ and $n_{ik}$ is the outer unit normal to
the $i$th side of $B_k$.

It remains to show that $\nu_k$ weakly converges to $\nu$. For any
continuous function $f$ on $\TTT^3$ one has
$$
\int\!\!\!\!\int\!\!\!\!\int_{\TTT^3} f(v, v^+, n)\, d\nu_k(v, v^+,
n)\ =\ \sum_i |b_k^i| \int\!\!\!\!\int_{\TTT^2} f(v, v^+, n_{ik})\,
d\eta_{ik}(v, v^+)\ =\
$$
 \beq\lab{sum fraction}
=\ \sum_i \frac{|b_k^i|}{|\pl B_k^i|}
\int\!\!\!\!\int\!\!\!\!\int_{\TTT^2 \times \SSS_k^i}
f(\Upsilon_{ik}(v, v^+, n), n_{ik})\, d\nu(v, v^+, n).
 \eeq
For each $k$ define the mapping from $\TTT^3$ to $\TTT^3$ by the
relations~ $(v, v^+, n) \mapsto (\Upsilon_{ik}(v, v^+, n), n_{ik})$
if $n \in \SSS_k^i$. It uniformly converges to the identity mapping
as $k \to \infty$; hence the function $\tilde f_k$, defined by the
relations $\tilde f_k(v, v^+, n) := f(\Upsilon_{ik}(v, v^+, n),
n_{ik})$ if $n \in \SSS_k^i\,$, uniformly converges to $f$ as $k \to
\infty$. From here and from (\ref{max i}) it follows that the right
hand side in (\ref{sum fraction}) converges to
$\int\!\!\!\int\!\!\!\int_{\TTT^3} f(v, v^+, n)\, d\nu(v, v^+, n)$
as $k \to \infty$. Thus, the convergence $\int f d\nu_k \to \int f
d\nu$ is proved. Q.E.D.
 \end{proof}

The inverse statement of the theorem follows from statements 1 and
2. Indeed, let $\nu \in \MMM_{B}$. Using statement \ref{p_2}, find a
sequence of convex polygons $B_k \subset B$ and a sequence $\nu_k
\in \MMM_{B_k}$ weakly converging to $\nu$. According to statement
\ref{p_1}, each measure $\nu_k$ is generated by a rough body.
Consider the sequence of sets $\QQQ_{kl} \subset B_k$,\, $l = 1,\,
2, \ldots$ representing this body, and then from all these sequences
choose a diagonal sequence $\tilde \QQQ_k = \QQQ_{kl_k}$ such that
the corresponding sequence of measures $\nu_{\tilde\QQQ_k, B}$
weakly converges to $\nu$ and Area$(B \setminus \tilde\QQQ_{k})$
goes to zero as $k \to \infty$. The sequence $\tilde \QQQ_k$
represents a body $\BBB$ obtained by roughening $B$ and $\nu_\BBB =
\nu$.

\section{Proof of the lemmas}

\subsection{Proof of lemma \ref{l1}}

Consider the billiard in $\RRR^2 \setminus \QQQ$. For $(\xi, v) \in
\set$, denote by $\tau(\xi, v)$ the time the billiard trajectory
with the initial data $\xi,\, -v$ spends in $B \setminus Q$. In
particular, if $\xi \in \pl B \cap \pl Q$, one has $\tau(\xi, v) =
0$.

Denote by $D$ the set of points $(x,w) \in (B \setminus Q) \times
S^1$ that are accessible from $\set$; that is, there exists $(\xi,
v) \in \set$ such that the billiard particle with the data $\xi,\,
-v$ at the zero moment of time, at some moment $0 \le t \le
\tau(\xi, v)$ will pass through $x$ with the velocity $w$. This
description defines the change of coordinates in $D:\ (\xi, v, t)
\mapsto (x, w);\ (\xi, v) \in \set,\ t \in [0,\, \tau(\xi, v)]$, and
the element of phase volume $d^2x\, dw$ in the new coordinates takes
the form $d\mu(\xi, v)\, dt$. Hence, the phase volume of $D$ equals
$\int\!\!\!\int\!\!\!\int_{D} d^2 x\, dw = \int\!\!\!\int_{\set}
\tau(\xi, v)\, d\mu(\xi, v)$. Taking into account that $D \subset (B
\setminus Q) \times S^1$ and the phase volume of $(B \setminus Q)
\times S^1$ equals $2\pi \cdot \text{Area}(B \setminus Q)$, one gets
\begin{equation}\label{ast1}
\int\!\!\!\int_{\set} \tau(\xi, v)\, d\mu(\xi,v) \le 2\pi \cdot
\text{Area}(B \setminus Q).
\end{equation}
This is in fact a simple modification of the well-known {\it mean
free path} formula (see, e.g., \cite{chernov}).

One has $\tau(\xi, v) \ge |\xi - \xi^+_{Q,B}(\xi, v)|$: the time the
particle spends in $B \setminus Q$ exceeds the distance between the
initial and final points of the trajectory. This inequality and
(\ref{ast1}) imply (a).

The points $\xi$ and $\xi^+_{Q,B}(\xi, v)$ divide the curve $\pl B$
into two arcs; denote by $\gam(\xi, v)$ the shortest one. One has
$|\gam(\xi, v)| = l(\xi, \xi^+_{Q,B}(\xi, v))$, therefore $|\xi -
\xi^+_{Q,B}(\xi, v)| \ge c \, |\gam(\xi, v)|$. It follows that
\begin{equation}\label{ast2}
c \int\!\!\!\!\!\!\!\!\!\!\int\limits_{\set} |\gam(\xi, v)|\,
d\mu(\xi,v) \le \int\!\!\!\!\!\!\!\!\!\!\int\limits_{\set} |\xi -
\xi^+_{Q,B}(\xi, v)|\, d\mu(\xi,v) \le 2\pi \cdot \text{Area}(B
\setminus Q).
\end{equation}

Let $\varrho(y)$ be a natural parametrization of the curve $\pl B,\
\varrho: [0,\, |\pl B|] \to \pl B$. By $f(y)$ denote the measure of
the values $(\xi, v)$ such that the interval $\gam(\xi, v)$ contains
the point $\varrho(y)$; that is, $f(y) := \int\!\!\!\int_{\set}
\mathbb{I} (\varrho(y) \in \gam(\xi, v))\, d\mu(\xi,v)$. Making
change of variables in the integral in the left hand side of
(\ref{ast2}), one gets
$$
\int\!\!\!\int_{\set} |\gam(\xi, v)|\, d\mu(\xi,v) = \int_0^{|\pl
B|} f(y)\, dy,
$$
therefore
\begin{equation}\label{ast3}
\int_0^{|\pl B|} f(y)\, dy \le \frac{2\pi}{c}\, \text{Area}(B
\setminus Q).
\end{equation}

One easily sees that $|f(y_1) - f(y_2)| \le 4\, |y_1 - y_2|$ for any
$y_1$ and $y_2$ and $f(y) \ge 0$. From here and from (\ref{ast3}) it
follows that for sufficiently small $\text{Area}(B \setminus Q)$
(namely, for $\text{Area}(B \setminus Q) \le {c} |\pl B|^2/(2\pi)$)
holds $f(y) \le \sqrt{8\pi/c}\, \sqrt{\text{Area}(B \setminus Q)}$.

Recall that $\Arg(v)$ is the angle the vector $v \ne 0$ forms with a
fixed vector $v_0$; the angle is measured clockwise from $v_0$ to
$v$ and is defined modulo\,$2\pi$. Introduce the shorthand notation
$\xi^+ := \xi^+_{Q,B}(\xi, v)$ and denote by $\Del\Arg(\xi, v)$ the
smallest in modulus of the values $\Arg(n(\xi^+)) - \Arg(n(\xi))$.
In other words, $\Del\Arg(\xi,v)$ equals to the smallest of the
values
$$
\int_{\gam(\xi,v)} |d\Arg(n_{\xi'})|,~~~~~ \int_{\pl B \setminus
\gam(\xi,v)} |d\Arg(n_{\xi'})|.
$$
Taking into account that $|n(\xi^+) - n(\xi)| \le |\Del\Arg(\xi,
v)|$, one gets that
$$
|n(\xi^+) - n(\xi)| \le \int_{\gam(\xi,v)} |d\Arg(n_{\xi'})|,
$$
and therefore,
$$
\overline{|n - n^+|}_{Q,B}\ \le\ \int\!\!\!\int_{\set} \left(
\int_{\gam(\xi,v)} |d\Arg(n_{\xi'})| \right) d\mu(\xi,v).
$$
Making change of variables in this integral, one obtains
$$
\overline{|n - n^+|}_{Q,B}\ \le\ \int_0^{|\pl B|} f(y)\
|d\,\Arg(n_{\varrho(y)})| \le 2\pi\, \sqrt{{8\pi}/{c}}\,
\sqrt{\text{Area}(B \setminus Q)}.
$$
Thus, (b) is also proved.

\subsection{Proof of lemma \ref{l2}}

Fix $n \in S^1$ and $m \in \mathbb N$. Let $\s$ be an involutive
permutation of $\{ 1, \ldots, m \}$, that is, $\s^2 =$ id. Divide
the half-circumference $S_n^1 := \{ v \in S^1 : \langle v,n \rangle
\ge 0 \}$ into $m$ arcs $\SSS_{n,m}^1 = \SSS_n^1$, \ldots,
$\SSS_{n,m}^m = \SSS_n^m$ numbered clockwise, such that for any
$i$,\, $\int_{\SSS_n^i} \langle v,n \rangle\, dv = 2/m$. For the
sake of brevity we omit the subscript $m$ when no confusion can
arise.

\begin{opr}\label{o5}
A measure $\eta$ is called a $(\s, n)$-measure if $\eta \in \Lam_n$
and spt$\,\eta \subset \cup_{i=1}^m \left( \SSS_n^i \times
\SSS_n^{\s(i)} \right)$, and therefore, for any $i$ holds $\eta
\left( \SSS_n^i \times \SSS_n^{\s(i)} \right) = 2/m$.
\end{opr}

\begin{utv}\label{u2}
For any measure $\eta \in \Lam_n$ there exists a sequence of
involutive permutations $\s_k$ on $\{ 1, \ldots, m_k \}$,\, $k =
1,\, 2, \ldots$ such that $m_k$ tends to infinity and any sequence
of $(\s_k, n)$-measures weakly converges to $\eta$ as $k \to
\infty$.
\end{utv}

\begin{utv}\label{u1}
Let $\s$ be an involutive permutation on $\{ 1, \ldots, m \}$. Then
the distance (in variation) between the set of measures generated by
$n$-hollows and the set of $(\s, n)$-measures does not exceed
$16/m$. In other words, whatever $\ve > 0$, there exist a $(\s,
n)$-measure $\eta$ and an $n$-hollow $(\Om, I)$ such that $\|
\eta_{\Om,I} - \eta \| < 16/m + \ve$; here the norm means variation
of measure.
\end{utv}

This distance actually equals zero, but we only need the (weaker)
claim of proposition \ref{u1}.

Lemma \ref{l2} follows from propositions \ref{u2} and \ref{u1}.
Indeed, let $\eta \in \Lam_n$. First choose the sequence of
permutations $\s_k$, according to proposition \ref{u2}, and then,
using proposition \ref{u1}, for every $k$ choose an an $n$-hollow
$(\Om_{k},I_{k})$ such that the distance from $\eta_{\Om_{k},I_{k}}$
to the set of $(\s_k, n)$-measures does not exceed $17/{m_k}$. The
sequence of chosen measures $\eta_{\Om_{k},I_{k}}$ weakly converges
to $\eta$.

\subsection{Proof of proposition \ref{u2}}

Introduce on $S^1_n$ the angular coordinate $\vphi = \text{Arg}\, v
- \text{Arg}\, n$; that is, $\vphi$ changes between $-\pi/2$ and
$\pi/2$ and increases clockwise. With this notation, to the arcs
$\SSS^i_{n,m}$ correspond the segments $J^i_m =
[\arcsin(-1+2(i-1)/m),\, \arcsin(-1+2i/m)]$. Define the measure
$\lam$ on $\interval$ by $d\lam(\vphi) = \cos\vphi\, d\vphi$ and
denote by $\Lam$ the set of measures $\eta$ on $\Box := \interval
\times \interval$ such that~ (a) $\pi_\vphi^\# \eta = \lam =
\pi_{\vphi^+}^\# \eta$ ~ and (b) $\pi_d^\# \eta = \eta$.~ Here
$\pi_\vphi$,\, $\pi_{\vphi^+}$, and $\pi_d$ are defined by
$\pi_\vphi(\vphi, \vphi^+) = \vphi$,\, $\pi_{\vphi^+}(\vphi,
\vphi^+) = \vphi^+$,\, $\pi_d(\vphi, \vphi^+) = (\vphi^+, \vphi)$.
Reformulating definition \ref{o5}, we shall say that $\eta$ is a
$\s$-measure if $\eta \in \Lam$ and spt$\,\eta \subset \cup_{i=1}^m
\left( J_m^i \times J_m^{\s(i)} \right)$. Notice that in the new
notation, the objects do not depend on $n$ anymore: we write $\Lam$
instead of $\Lam_n,\ \s$-measure instead of $(\s, n)$-measure, and
hollow in place of $n$-hollow.

In this notation, proposition \ref{u2} can be reformulated as
follows:~ for any measure $\eta \in \Lam$ there exists a sequence of
involutive permutations $\s_k$ on $\{ 1, \ldots, m_k \}$,\, $k =
1,\, 2, \ldots$ such that $m_k$ tends to infinity and any sequence
of $\s_k$-measures weakly converges to $\eta$ as $k \to \infty$.

The idea of the proof is as follows. First, $\eta$ is approximated
by means of a rational matrix, and then, this matrix is approximated
by means of a larger matrix generated by a permutation.

Consider the partition of $\Box$ into smaller rectangles
$\Box_k^{ij} = J^i_k \times J^j_k$,\, $i,\, j = 1, \ldots, k$.
Choose rational nonnegative numbers $c^{ij}_k$ such that $c^{ij}_k =
c^{ji}_k$, $\sum_j c^{ij}_k = 2/k$ for any $i$, and $\big|
\eta\left( \Box^{ij}_k \right) - c_k^{ij} \big| \le k^{-3}$ for any
$i$ and $j$. To do so, it suffices to take positive rational values
$c_k^{ij}$ such that $\eta(\Box_k^{ij}) - k^{-4} \le c_k^{ij} \le
\eta(\Box_k^{ij})$ for $i > j$ and put $c_k^{ij} = c_k^{ji}$ for $i
< j$ and $c_k^{ii} = 2/k - \sum_{j\ne i} c_k^{ij}$ for $i = j$. One
has $\eta\left(J_k^i \times \interval\right) = \sum_{j=1}^k
\eta(\Box_k^{ij}) = 2/k$, hence $c_k^{ii} - \eta(\Box_k^{ii}) =
\sum_{j\ne i} \left( \eta(\Box_k^{ij}) - c_k^{ij} \right) \in [0,\,
(k - 1) \cdot k^{-4}] \subset [0,\, k^{-3}]$.

Any sequence of measures $\eta_k$ satisfying the conditions
$\eta_k(\Box_k^{ij}) = c_k^{ij}$,\,\, $1 \le i,\, j \le k$ weakly
converges to $\eta$. Indeed, for any continuous function $f$ on
$\Box$ holds
$$
\int_\Box f\, d\eta_k - \int_\Box f\, d\eta = \sum_{i,j=1}^k
\int_{\Box_k^{ij}} f\, \left( d\eta_k - d\eta \right) \le k^{-1}
\max f \to 0
$$
as $k \to \infty$.

To complete the proof, it suffices to find an integer $m_k > k$ and
an involutive permutation $\s_k$ of $\{ 1, \ldots, m_k \}$ such that
any $\s_k$-measure, $\eta_k$, satisfies the equalities
$\eta_k(\Box_k^{ij}) = c_k^{ij}$,\, $i,\, j = 1, \ldots, k$. Choose
a positive integer $N$ such that all the values $a_{ij} := N \cdot
c_k^{ij}$ are integer. The obtained matrix $A = (a_{ij})_{i,j=1}^k$
is symmetric and for any $i$ the value $\sum_{j=1}^k a_{ij} = 2N/k$
is a fixed positive integer. In Appendix B it is shown that there
exist square matrices $B_{ij} = (b_{ij}^{\mu\nu})_{\mu,\nu}$ of size
$2N/k$ such that $B_{ij}^T = B_{ji}$, the sum of elements in any
matrix $B_{ij}$ equals $a_{ij}$ and the block matrix $D = (B_{ij})$
composed of these matrices has exactly one unit in each row and each
column, and other elements are zeros.

$D$ is a symmetric square matrix of size $2N$; denote its elements
by $d_{ij}$. Put $m_k = 2N$ and define the mapping $\s_k$ on $\{ 1,
\ldots, 2N \}$ in such a way that $d_{i\s_k(i)} = 1$ for any $i$.
The so defined mapping $\s_k$ is a permutation; it is involutive
since the matrix $D$ is symmetric. Moreover, if $\eta_k$ is a
$\s_k$-measure then for any $i$ and $j$ holds $\eta_k(\Box_k^{ij}) =
N^{-1} \sum_{\mu,\nu} b_{ij}^{\mu\nu} = c_k^{ij}$. The proposition
is proved.

\subsection{Proof of proposition \ref{u1}}

{\bf 1.}~ Whatever the $n$-hollow $(\Om, I)$, one introduces the
reference system $(x_1, x_2)$ in such a way that $n$ coincides with
$(0, -1)$, and the interval $I$ belongs to the straight line $x_2 =
0$ and contains the origin $O = (0,0)$. Like in the proof of
proposition \ref{u2}, introduce the coordinate $\vphi = \text{Arg}\,
v - \text{Arg}\, n$ on $S^1_n$. One has $v = -(\sin\vphi,
\cos\vphi)$,\, $\vphi \in [-\pi/2,\, \pi/2]$. The definition of the
segments $J_m^i = J^i$, the measure $\lam$, the set of measures
$\Lam$, and the $\s$-measure see in the beginning of the previous
subsection. The mapping $(\xi, v) \mapsto V_{\Om,I}^+(\xi, v)$ in
the new coordinates $\xi,\, \vphi$ is written as $(\xi, \vphi)
\mapsto \vphi_{\Om,I}^+(\xi, \vphi)$. Finally, define the measure
$\mu_I$ on $I \times \interval$ by $d\mu_I(\xi, \vphi) =
\frac{\cos\vphi}{|I|}\, d\xi\, d\vphi$.

Denote $\Box^{\,'} = \left( \cup_{i=2}^{m-1} J^i \right) \times
\left( \cup_{i=2}^{m-1} J^i \right),\ \Box_1 = J^1 \times
\interval$,\, $\Box_2 = J^m \times \interval$,\, $\Box_3 = \left(
\cup_{i=2}^{m-1} J^i \right) \times J^1$, and $\Box_4 = \left(
\cup_{i=2}^{m-1} J^i \right) \times J^m$. Thus, one has $\Box
\setminus \Box^{\,'} = \Box_1 \cup \Box_2 \cup \Box_3 \cup \Box_4$;
see the figure below.
 
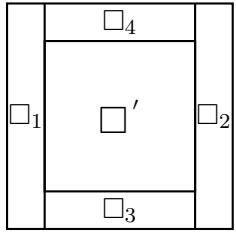
\vspace*{5mm}

 \rput(5,-1.5){
 \pspolygon(-1,-1)(1,-1)(1,1)(-1,1)
 \pspolygon(-1.5,-1.5)(1.5,-1.5)(1.5,1.5)(-1.5,1.5)
 \psline(-1,-1.5)(-1,1.5)
 \psline(1,-1.5)(1,1.5)
 \rput(-1.25,0){\small $\Box_1$}
 \rput(1.25,0){\small $\Box_2$}
 \rput(0,-1.25){\small $\Box_3$}
 \rput(0,1.25){\small $\Box_4$}
 \rput(0,0){\large $\Box^{\,'}$}

 }

 \vspace{35mm}

It suffices to construct a sequence of hollows $(\Om_\ve, I_\ve)$,\,
$\ve > 0$ such that
 \begin{quote}
{\bf (P)}~ {\it for any $i \ne 1,\ m,\ \s(1),\ \s(m)$ the measure of
the set of values $(\xi, \vphi) \in I_\ve \times J^i$ such that
$\vphi^+_{\Om_\ve,I_\ve} (\xi, \vphi) \not\in J^{\s(i)}$ goes to
zero as $\ve \to 0$}.
 \end{quote}
Then, speaking of restrictions of measures on the subset
$\Box^{\,'}$, one gets that the distance from the restrictions of
measures $\eta_{\Om_\ve,I_\ve}$ to the set of restrictions of
$\s$-measures goes to zero as $\ve \to 0$. On the other hand, for
any measure $\eta \in \Lam$ one has $\eta(\Box_1) = \eta(\Box_2) =
2/m$,\, $\eta(\Box_3) \le 2/m$,\ $\eta(\Box_4) \le 2/m$, hence
$\eta(\Box \setminus \Box^{\,'}) \le 8/m$, therefore the distance
between the restrictions on $\Box \setminus \Box^{\,'}$ of any two
measures $\eta_1$ and $\eta_2$ from $\Lam$ does not exceed $16/m$:
$\ \|\, \eta_1 \rfloor_{\Box \setminus \Box^{\,'}} -
\eta_2\rfloor_{\Box \setminus \Box^{\,'}}\, \| \le 16/m$. It follows
that the upper limit of distances from $\eta_{\Om_\ve,I_\ve}$ to the
set of $\s$-measures does not exceed $16/m$, and so, proposition
\ref{u1} is proved.

{\bf 2.}~ The rest of this subsection is dedicated to the detailed
description of the sequence of hollows $(\Om_\ve, I_\ve)$ and to the
proof of property (P) for them.

First consider an auxiliary construction. Take two different points
$F$ and $F\marka$ above the line $l = \{ x_2 = 0 \}$, with $|OF| = 2
= |OF\marka|$. Denote by $\Phi$ and $\Phi\marka$ the angles the rays
$OF$ and $OF\marka$, respectively, form with the vector $(0,1)$. The
angles are counted clockwise from $(0,1)$. Thus, $F = 2(\sin\Phi,
\cos\Phi)$ и $F\marka = 2(\sin\Phi\marka, \cos\Phi\marka)$. Assume,
for further convenience, that $F$ is situated on the left of
$F\marka$; thus, one has $-\pi/2 < \Phi < \Phi\marka < \pi/2$. (The
case where $F$ is situated on the right of $F\marka$ is completely
similar.) Select three positive numbers $\lam$,\, $\lam\marka$, and
$\del$, and define two ellipses $\EEE$ and $\EEE\marka$ and two
parabolas $\PPP$ and $\PPP\marka$. The first ellipse has the foci
$O$ and $F$, the length of its large semiaxis is $\sqrt{1 + \lam}$,
of the small semiaxis, $\sqrt\lam$, and the focal distance equals 2.
The second ellipse has the foci $O$ and $F\marka$, the lengths of
its large and small semiaxes are $\sqrt{1 + \lam\marka}$ and
$\sqrt{\lam\marka}$, respectively, and the focal distance is also 2.
The parabolas $\PPP$ and $\PPP\marka$ have the foci $F$ and
$F\marka$, respectively, the common axis $FF\marka$, and the same
focal distance $\del$. Thus, the parabolas are symmetric to each
other with respect to the bisectrix of the triangle $OFF\marka$. The
parameter $\del$ is chosen sufficiently small, so that the point $O$
lies in the exterior of both parabolas. See the figure below.

 
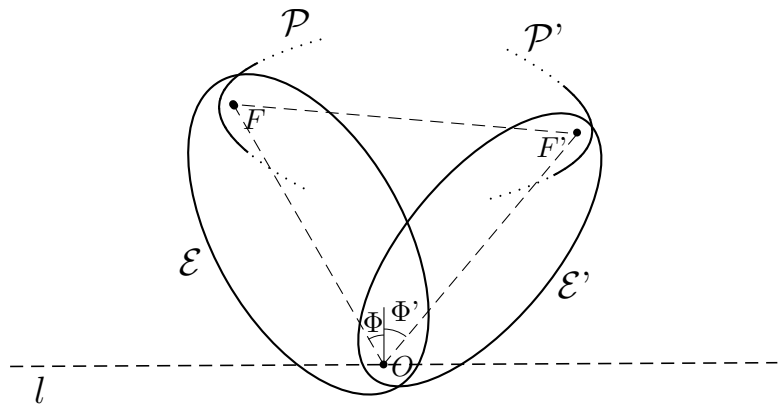
\vspace{0mm}

 \vspace*{38mm}

 \rput{120}(5.8,0.2){
 \psellipse(0,0)(2.366,1.264)
 \psdots[dotsize=3pt](2,0)(-2,0)
  \psline[linewidth=.4pt,linestyle=dashed](2,0)(-2,0)
   }
 \rput{50}(8,0){
 \psellipse(0,0)(2.236,1)
  \psline[linewidth=.4pt,linestyle=dashed](2,0)(-2,0)
 }
 \psline[linewidth=.6pt,linestyle=dashed](1.8,-1.5)(12,-1.45)

        \rput{-6.2}(8,2.2){
 \rput{-90}(-3.5,0){
 \parabola(0.6,0.45)(0,0)
 \parabola[linestyle=dotted](1,1.25)(0,0)
 \psdot[dotsize=3pt](0,0.2)
 }
 \rput{90}(1.47,0){
 \parabola(0.6,0.45)(0,0)
 \parabola[linestyle=dotted](1,1.25)(0,0)
 \psdot[dotsize=3pt](0,0.2)
 }
  \psline[linewidth=.4pt,linestyle=dashed](-3.2,-0.05)(1.15,0.05)
  }
  \rput(2.2,-1.3){\large $l$}
  \rput(4.2,.4){\large $\EEE$}
  \rput(9.3,.1){\large $\EEE\marka$}
 \rput(5.6,3.6){\large $\PPP$}
 \rput(8.9,3.4){\large $\PPP\marka$}
 \rput(5.05,2.3){\small $F$}
 \rput(9,1.93){\small $F\marka$}
 \rput(7.02,-0.98){\small $O$}
   \psline[linewidth=0.4pt](6.77,-0.97)(6.77,-0.2)
 \psarc[linewidth=0.2pt](6.77,-0.97){0.4}{90}{121}
 \psarc[linewidth=0.2pt](6.77,-0.97){0.48}{52}{90}
  \rput(6.63,-0.4){\small $\Phi$}
  \rput(7.05,-0.25){\small $\Phi\marka$}

 \vspace{28mm}

In what follows, we shall distinguish between the billiard and {\it
pseudo-billiard} dynamics. The pseudo-billiard dynamics is defined
as follows. A particle starts at a point $(\xi, 0) \in l$ and moves
with a velocity $(\sin\vphi, \cos\vphi)$ until it reflects from the
interior side of $\EEE$. (Before the reflection it can intersect
other curves $\EEE\marka$,\, $\PPP$,\, $\PPP\marka$, or even
intersect $\EEE$ from the outer side, without changing the
velocity.) Then it moves again with constant velocity, until it
reflects from the interior side of $\PPP$. Then, in the same way, it
reflects from the interior side of $\PPP\marka$ and then from the
interior side of $\EEE\marka$, and finally, intersects $l$ from
above to below. Denote by $(\xi\marka, 0)$ the point of
intersection, and by $-(\sin\vphi\marka, \cos\vphi\marka)$, the
velocity at this point.

Consider the admissible set: the set of 7-uples $(\vphi, \xi, \Phi,
\Phi\marka, \lam, \lam\marka, \del)$ such that all the indicated
reflections occur in the prescribed order. This set is open and
nonempty. Indeed, let $\del(\Phi, \Phi\marka)$ be the least of the
values $\del$ such that one of the parabolas (in fact, both of them
simultaneously) passes through $O$. Put $\vphi = \Phi$,\, $\xi = 0$,
and take arbitrary values $\lam > 0$,\, $\lam\marka > 0$,\, $-\pi/2
< \Phi < \Phi\marka < \pi/2$,\, $0 < \del < \del(\Phi, \Phi\marka)$.
The particle with initial data $\vphi = \Phi$,\, $\xi = 0$ first
passes along the large semiaxis of $\EEE$, then reflects from
$\EEE$, returns along the same semiaxis and reflects from $\PPP$.
Then it moves with the velocity parallel to $FF\marka$, reflects
from $\PPP\marka$, moves the large semiaxis of the ellipse
$\EEE\marka$, reflects from it and returns to $O$ along the same
semiaxis. Thus, the admissible set is nonempty. Under a small
perturbation of the parameters $\vphi$,\, $\xi$,\, $\Phi$,\,
$\Phi\marka$,\, $\lam$,\, $\lam\marka$,\, $\del$, all the
reflections are maintained and the order of reflections remains the
same. This implies that the admissible set is open.

This description determines the mapping $\vphi\marka =
\vphi\marka(\vphi, \xi, \Phi, \Phi\marka, \lam, \lam\marka,
\del)$,\, $\xi\marka = \xi\marka(\vphi, \xi, \Phi, \Phi\marka, \lam,
\lam\marka, \del)$\footnote{Note that throughout this paper the sign
{\marka} (prime) never means derivation.} from the admissible set to
$\RRR^2$. This mapping is infinitely differentiable. For $\vphi =
\Phi$ and $\xi = 0$ one has
 \beq\lab{1 L2}
\vphi\marka(\Phi, 0, \Phi, \Phi\marka, \lam, \lam\marka, \del) =
\Phi\marka.
 \eeq
For $\xi = 0$ with arbitrary $\vphi$ one has
 \beq\lab{2 L2}
\xi\marka(\vphi, 0, \Phi, \Phi\marka, \lam, \lam\marka, \del) = 0,
 \eeq
and
$$
\vphi\marka(\vphi, 0, \Phi, \Phi\marka, \lam, \lam\marka, \del)~~
\text{ does not depend on }~ \del.
$$

Indeed, a particle starting at $O$, after the reflection from $\EEE$
passes through $F$, then after reflecting from $\PPP$ moves in
parallel to $FF\marka$, after the reflection from $\PPP\marka$
passes through $F\marka$, and finally, after the reflection from
$\EEE\marka$ returns to $O$ (see the figure below). The initial and
final velocity of the particle are, respectively, $(\sin\vphi,
\cos\vphi)$ and $-(\sin\vphi\marka, \cos\vphi\marka)$. Denoting by
$\al$ and $\al\marka$ the angles the second and fourth segments of
the (5-segment) trajectory form, respectively, with $OF$ and
$OF\marka$, one has $\al = \al\marka$. The angle $\al$ is a function
of $\vphi$, and $\vphi\marka$ is a function of $\al\marka$; these
functions depend only on the parameters of the ellipses $\EEE$ and
$\EEE\marka$, respectively, and do not depend on the parameter
$\del$ determining the shape of parabolas.
 
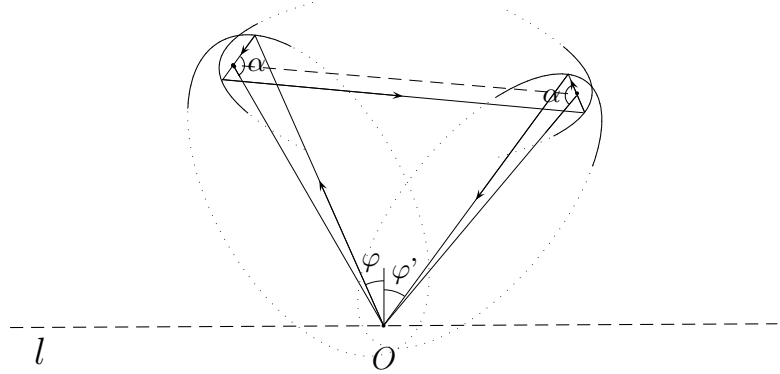
\vspace{35mm}

 \rput{120}(5.8,0.2){
 \psellipse[linewidth=.4pt](0,0)(2.366,1.264)
 \psdots[dotsize=1.5pt](2,0)(-2,0)
 \psline[linewidth=.4pt](2,0)(-2,0)
 }
 \rput{50}(8,0){
 \psellipse[linewidth=.4pt](0,0)(2.236,1)
  \psline[linewidth=.2pt](2,0)(-2,0)
  }
    \pspolygon[fillstyle=solid,linecolor=white,fillcolor=white](4,-2)(9.5,-2)(10,0.3)(6.5,2.8)(4,1.3)
 \psline[linewidth=.4pt,linestyle=dashed](1.8,-1.5)(12,-1.45)
   \rput{120}(5.8,0.2){
 \psellipse[linewidth=.4pt,linestyle=dotted](0,0)(2.366,1.264)
 \psdots[dotsize=1.5pt](2,0)(-2,0)
 \psline[linewidth=.15pt](2,0)(-2,0)
 \psline[linewidth=.2pt](1.9,0.225)(2.2,-0.45)(-2,0)
 \psline[linewidth=.2pt,arrows=<-,arrowscale=1.5](0.1,-0.225)(-2,0)
 \psline[linewidth=.2pt,arrows=<-,arrowscale=1.5](2.067,-0.15)(2.2,-0.45)
    \psarc[linewidth=0.4pt](2,0){0.15}{177}{290}
      }
 \rput{50}(8,0){
 \psellipse[linestyle=dotted,linewidth=.4pt](0,0)(2.236,1)
  \psline[linewidth=.15pt](2,0)(-2,0)
 \psline[linewidth=.2pt](1.9,-0.225)(2.15,0.27)(-2,0)
 \psline[linewidth=.2pt,arrows=->,arrowscale=1.5](1.9,-0.225)(2.1143,0.1993)
  \psline[linewidth=.2pt,arrows=->,arrowscale=1.5](2.15,0.27)(0.075,0.135)
   \psarc[linewidth=0.4pt](2,0){0.15}{67}{177}
   }

        \rput{-6.2}(8,2.2){
 \rput{-90}(-3.5,0){
 \parabola[linewidth=.4pt](0.6,0.45)(0,0)
 \parabola[linestyle=dotted,linewidth=.4pt](1,1.25)(0,0)
 \psdot[dotsize=1.5pt](0,0.2)
 }
 \rput{90}(1.47,0){
 \parabola[linewidth=.4pt](0.6,0.45)(0,0)
 \parabola[linestyle=dotted,linewidth=.4pt](1,1.25)(0,0)
 \psdot[dotsize=1.5pt](0,0.2)
  }
  \psline[linewidth=.2pt,linestyle=dashed](-3.2,-0.05)(1.15,0.05)
  \psline[linewidth=.2pt](-3.44,-0.26)(1.385,-0.18)
  \psline[linewidth=.2pt,arrows=->,arrowscale=1.5](-3.44,-0.26)(-1.0275,-0.22)
    }
  \rput(2.2,-1.3){\large $l$}
  \rput(6.77,-1.4){$O$}
 \psline[linewidth=0.2pt](6.77,-0.97)(6.77,-0.2)
 \psarc[linewidth=0.2pt](6.77,-0.97){0.6}{90}{114}
 \psarc[linewidth=0.2pt](6.77,-0.97){0.48}{55}{90}
  \rput(6.6,-0.13){\small $\vphi$}
  \rput(7.05,-0.23){\small $\vphi\marka$}
     \rput(9,2.1){\small $\al$}
     \rput(5.08,2.5){\small $\al$}

 \vspace{27mm}

Using properties of ellipses, one derives the formulas connecting
$\vphi$,\, $\al$, and $\vphi\marka = \vphi\marka(\vphi, 0, \Phi,
\Phi\marka, \lam, \lam\marka, \del)$:
 \beq\lab{3 L2}
\sin(\vphi-\Phi) = \frac{\lam\sin\al}{2 + \lam - 2\cos\al \sqrt{1 +
\lam}}\,,~~~ \sin(\vphi\marka - \Phi\marka) =
-\frac{\lam\marka\sin\al}{2 + \lam\marka - 2\cos\al \sqrt{1 +
\lam\marka}}\,.
 \eeq
It follows that
 \beq\lab{4 L2}
\frac{\pl\vphi\marka}{\pl\vphi}\Bigg\rfloor_{\!\!\!\scriptsize
\begin{array}{l} \vphi=\Phi\\ \xi = 0
\end{array}} = -\left( \frac{\sqrt{\lam\marka}}{1 + \sqrt{\lam\marka}}\ \frac{1 +
\sqrt{\lam}}{\sqrt{\lam}} \right)^2.
 \eeq
With fixed $\Phi$,\, $\Phi\marka$,\, $\lam$,\, $\lam\marka$, and
$\del$ the mapping $\vphi\marka(\vphi, \xi)$,\, $\xi\marka(\vphi,
\xi)$ preserves the measure, $\cos\vphi\, d\vphi\, d\xi =
\cos\vphi\marka\, d\vphi\marka\, d\xi\marka$, hence
$$
\cos\vphi\, =\, \pm \cos\vphi\marka\, \left| \begin{array}{cc}
\frac{\pl\vphi\marka}{\pl\vphi} & \frac{\pl\vphi\marka}{\pl\xi} \\
\frac{\pl\xi\marka}{\pl\vphi} & \frac{\pl\xi\marka}{\pl\xi}
\end{array} \right|\,.
$$
Using (\ref{2 L2}), one gets that
$\frac{\pl\xi\marka}{\pl\vphi}\big\rfloor_{\xi=0} = 0$, hence
$$
\left| \begin{array}{cc}
\frac{\pl\vphi\marka}{\pl\vphi} & \frac{\pl\vphi\marka}{\pl\xi} \\
\frac{\pl\xi\marka}{\pl\vphi} & \frac{\pl\xi\marka}{\pl\xi}
\end{array} \right|_{\xi=0} = \ \frac{\pl\vphi\marka}{\pl\vphi}\
\frac{\pl\xi\marka}{\pl\xi}\bigg\rfloor_{\xi=0}\ ,
$$
therefore
 \beq\lab{5 L2}
\cos\vphi\, =\, \pm \cos\vphi\marka\
\frac{\pl\vphi\marka}{\pl\vphi}\
\frac{\pl\xi\marka}{\pl\xi}\bigg\rfloor_{\xi=0}\ .
 \eeq
Putting $\vphi = \Phi$,\, $\xi = 0$, and taking into account (\ref{1
L2}), (\ref{4 L2}) and (\ref{5 L2}), one gets
 \beq\lab{6 L2}
\cos\Phi\, =\, \pm \cos\Phi\marka  \left( \frac{\sqrt{\lam\marka}}{1
+ \sqrt{\lam\marka}}\ \frac{1 + \sqrt{\lam}}{\sqrt{\lam}} \right)^2
\frac{\pl\xi\marka}{\pl\xi}\Bigg\rfloor_{\!\!\!\scriptsize
\begin{array}{l} \vphi=\Phi\\ \xi = 0
\end{array}} .
 \eeq
Define the positive continuous functions $\lam(\Phi\marka)$ and
$\lam\marka(\Phi)$ by the relations
 \beq\lab{7 L2}
\left( \frac{\sqrt{\lam}}{1 + \sqrt{\lam}} \right)^{\!\!2} =\, \frac
12\, \cos\Phi\marka,~~~~ \left( \frac{\sqrt{\lam\marka}}{1 +
\sqrt{\lam\marka}} \right)^{\!\!2} =\, \frac 12\, \cos\Phi,
 \eeq
 then one has
 \beq\lab{8 L2}
\Bigg|\, \frac{\pl\xi\marka}{\pl\xi}\Bigg\rfloor_{\!\!\!\scriptsize
\begin{array}{c} \vphi=\Phi;\ \lam=\lam(\Phi\marka)\\ \xi=0;\ \lam\marka=\lam\marka(\Phi) \end{array}} \Bigg|\, =\, 1.
 \eeq
Additionally, taking into account (\ref{4 L2}) and (\ref{7 L2}), one
gets
 \beq\lab{9 L2}
\frac{\cos\Phi\marka}{\cos\Phi}\
\frac{\pl\vphi\marka}{\pl\vphi}\Bigg\rfloor_{\!\!\!\scriptsize
\begin{array}{c} \vphi=\Phi;\ \lam=\lam(\Phi\marka)\\ \xi=0;\ \lam\marka=\lam\marka(\Phi) \end{array}}\!\! =\, -1.
 \eeq
Recall that $\vphi\marka = \vphi\marka(\vphi, 0, \Phi, \Phi\marka,
\lam, \lam\marka)$, that is, the restriction of the function
$\vphi\marka$ to the subspace $\xi = 0$, does not depend on $\del$.
Hence the function
$\frac{\pl\vphi\marka}{\pl\vphi}\big\rfloor_{\xi=0}$ and, by formula
(\ref{5 L2}), the function
$\frac{\pl\xi\marka}{\pl\xi}\big\rfloor_{\xi=0}$ also do not depend
on $\del$. Put $\Phi_0 = \arcsin(1 - 2/m)$, so that $J^1 =
[-\pi/2,\, -\Phi_0]$ and $J^m = [\Phi_0,\, \pi/2]$, and put
$\Del\Phi = 2/m$. The set $\{ (\Phi, 0, \Phi, \Phi\marka,
\lam(\Phi\marka), \lam\marka(\Phi)):\ -\Phi_0 \le \Phi,\, \Phi\marka
\le \Phi_0,\ \Phi\marka - \Phi \ge \Del\Phi \}$ is compact and
belongs to the (open) domain of the function $\vphi\marka$. Choose a
sufficiently large integer value $k = k(\ve)$, so that for
 \beq\lab{51 L2}
 \begin{array}{cc}
 |\sin\vphi -
\sin\Phi| < {2}/(km),~~~~~  \xi = 0,~~~  -\Phi_0 \le \Phi,\,
\Phi\marka \le \Phi_0,\\
\Phi\marka - \Phi \ge \Del\Phi,~~~  \lam = \lam(\Phi\marka),~~~~~
\lam\marka = \lam\marka(\Phi)
 \end{array}
 \eeq
 holds true
 \beq\lab{10 L2}
-\frac{\cos\Phi\marka}{\cos\Phi}\ \frac{\pl\vphi\marka}{\pl\vphi}
\in [(1 + \ve)^{-1},\ 1 + \ve].
 \eeq
Formulas (\ref{10 L2}) and (\ref{1 L2}) mean that under the
conditions (\ref{51 L2}), $\vphi\marka$ is also close to
$\Phi\marka$. Increasing $k$ if necessary, ensure (under the same
conditions) that
 \beq\lab{11 L2}
\frac{\cos\Phi\marka}{\cos\Phi}\ \frac{\cos\vphi}{\cos\vphi\marka}
\in [(1 + \ve)^{-1},\ 1 + \ve].
 \eeq
Taking into account (\ref{5 L2}),\, (\ref{10 L2}), and (\ref{11
L2}), one obtains that under the conditions (\ref{51 L2}) holds true
 \beq\lab{12 L2}
\left|\, \frac{\pl\xi\marka}{\pl\xi}\, \right| \in [(1 + \ve)^{-2},\
(1 + \ve)^2].
 \eeq

\vspace{2mm}

{\bf 3.}~ Now we proceed to the description of the hollow $(\Om_\ve,
I_\ve)$.
  \vspace{1mm}

{\bf (a)}~  If $2 \le i \ne \s(i) \le m - 1$, divide the interval
$J^i$ into $k$ subintervals $J^{i,j}$ of equal measure $\lam$, going
in increasing order: $J^i = \cup_{j=1}^k J^{i,j}$,\, $\lam(J^{i,j})
= 2/(km)$ for any $j = 1, \ldots, k$. Recall that $d\lam(\vphi) =
\cos\vphi\, d\vphi$ and the value $k = k(\ve)$ is defined above.
Without loss of generality assume that $k(\ve) \to \infty$ as $\ve
\to 0$.

To each pair of intervals, $J^{i,j}$ and $J^{\s(i),j}$, we apply the
construction described above, see fig.\,9. Namely, draw arcs of
ellipses $\EEE_{i,j} = \duga{AB}$,\, $\EEE\marka_{\!\!i,j} =
\stackrel{\smile}{A\marka B\marka}$ and arcs of parabolas
$\PPP_{i,j}$,\, $\PPP\marka_{\!\!i,j}$. Without loss of generality
suppose that $i < \s(i)$. The angles $AOB$ and $A\marka OB\marka$
correspond to the angular intervals $J^{i,j}$ and $J^{\s(i),j}$,
respectively. The foci $\bar F = F_{i,j}$ and $\bar F\marka =
F\marka_{\!\!i,j}$ belong to the intervals $OA$ and $OA\marka$,
respectively. The endpoints of the arcs $\PPP_{i,j}$ and
$\PPP\marka_{\!\!i,j}$ also belong to the intervals $OA$ and
$OA\marka$, respectively. The angle corresponding to the ray $OA$
(and therefore to the left endpoint of the interval $J^{i,j}$) will
be denoted by $\bar\Phi = \Phi_{i,j}$, and the angle corresponding
to the ray $OA\marka$ (and therefore to the right endpoint of the
interval $J^{\s(i),j}$) will be denoted by $\bar\Phi\marka =
\Phi\marka_{\!\!i,j}$. Denote $\bar\lam = \lam_{i,j} :=
\lam(\bar\Phi\marka)$ and $\bar\lam\marka = \lam\marka_{i,j} :=
\lam\marka(\bar\Phi)$, according to the formula (\ref{7 L2}). Next,
select a value $\bar\del = \del_{i,j}$ and draw two curves ({\it
lateral reflectors}) in such a way that (i) each of the curves
contains an arc of parabola (the first curve contains $\PPP_{i,j}$,
and the second one, $\PPP\marka_{\!\!i,j}$), an arc of circumference
centered at $O$, and three radial segments;~ (ii) these curves do
not intersect the intervals whose endpoints belong to the set $\{
F_{\al, \bt},\ F\marka_{\!\! \gam,\del}:\ (\al, \bt) \neq (i, j),\,
(\gam, \del) \neq (\s(i), j) \}$: this will guarantee free passage
of particles from one parabola to another; and~ (iii) the
$\lam$-measure of the angular interval occupied by each lateral
reflector does not exceed $\ve/(km)$. On the figure below, the
angular reflectors are the curves joining the points $A$ and $C$,
and the points $A\marka$ and $C\marka$.

 
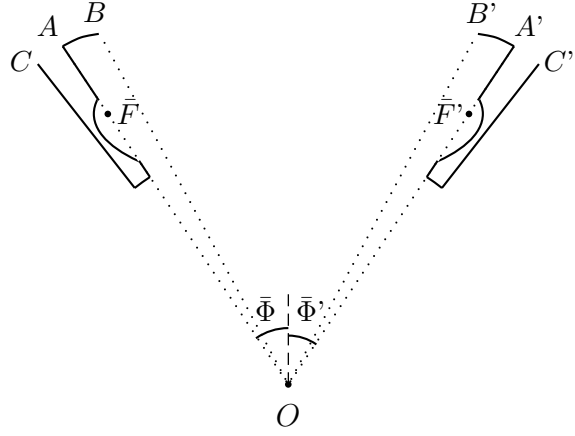
\vspace{0mm}

 \vspace*{55mm}

 \rput(6,0){

 \psellipse(0,3.6)(2.6,1)
 \pspolygon[fillstyle=solid,fillcolor=white,linecolor=white](0,0)(3,4.5)(0,5)(-3,4.5)

           \psdots[dotsize=2.5pt](0,0)(2.4,3.6)(-2.4,3.6)
   \pscurve(3,4.5)(2.76,4.616)(2.52,4.664)
   \psline[linestyle=dotted](3,4.5)(0,0)(2.5,4.65)
       \psline(1.84,2.76)(1.984,2.976)
       \psline(2.52,3.78)(3,4.5)
       \psarc(0,0){3.317}{52}{56.3}
       \psline(2.042,2.614)(3.331,4.263)
 \pscurve(-3,4.5)(-2.76,4.616)(-2.52,4.664)
   \psline[linestyle=dotted](-3,4.5)(0,0)(-2.5,4.65)
       \psline(-1.84,2.76)(-1.984,2.976)
       \psline(-2.52,3.78)(-3,4.5)
       \psarc(0,0){3.317}{123.7}{128}
       \psline(-2.042,2.614)(-3.331,4.263)
       \rput(0,-0.4){$O$}
       \rput(3.2,4.77){\small $A\marka$}
       \rput(3.6,4.3){\small $C\marka$}
       \rput(2.57,4.95){\small $B\marka$}
       \rput(2.13,3.63){\small $\bar F\marka$}
       \rput(-3.2,4.77){\small $A$}
       \rput(-3.55,4.3){\small $C$}
       \rput(-2.57,4.95){\small $B$}
       \rput(-2.13,3.63){\small $\bar F$}
 \psline[linestyle=dashed,linewidth=0.5pt](0,0)(0,1.2)
 \psarc[linewidth=0.8pt](0,0){0.75}{90}{125}
 \psarc[linewidth=0.8pt](0,0){0.65}{55}{90}
  \rput(-0.3,1){\small $\bar\Phi$}
  \rput(0.3,1){\small $\bar\Phi\marka$}

  }

 \vspace{15mm}

Notice that $-\Phi_0 \le \bar\Phi,\, \bar\Phi\marka \le \Phi_0$ and
$\bar\Phi\marka - \bar\Phi \ge \Del\Phi$. Indeed, $\bar\Phi$ and
$\bar\Phi\marka$ do not belong to the intervals $J^1 = [-\pi/2,\,
-\Phi_0]$ and $J^m = [\Phi_0,\, \pi/2]$. On the other hand, one has
$\bar\Phi\marka - \bar\Phi \ge \sin\bar\Phi\marka - \sin\bar\Phi =
\lam([\bar\Phi,\, \bar\Phi\marka]) \ge 2/m = \Del\Phi$.

Introduce the shorthand notation $\vphi\marka(\vphi, \xi) =
\vphi\marka(\vphi, \xi, \Phi_{i,j}, \Phi\marka_{i,j}, \lam_{i,j},
\lam\marka_{i,j}, \del_{i,j})$. According to (\ref{10 L2}) and
(\ref{12 L2}), for $\vphi \in J^{i,j}$ holds true
 \beq\lab{13 L2}
-\frac{\cos\bar\Phi\marka}{\cos\bar\Phi}\
\frac{\pl\vphi\marka}{\pl\vphi}(\vphi, 0) \in [(1 + \ve)^{-1},\ 1 +
\ve]
 \eeq
and
 \beq\lab{14 L2}
\left| \frac{\pl\xi\marka}{\pl\xi}(\vphi, 0) \right| \in [(1 +
\ve)^{-2},\ (1 + \ve)^2].
 \eeq
According to (\ref{1 L2}), one has $\vphi\marka(\bar\Phi, 0) =
\bar\Phi\marka$; this equality and the formula (\ref{13 L2}) imply
that for $\vphi \in J^{i,j}$ and $\vphi\marka = \vphi\marka(\vphi,
0)$ one has
 \beq\lab{15 L2}
-\frac{\cos\bar\Phi\marka}{\cos\bar\Phi}\ \frac{\vphi\marka -
\bar\Phi\marka}{\vphi - \bar\Phi} \in [(1 + \ve)^{-1},\ 1 + \ve].
 \eeq
On the other hand, one has
 \beq\lab{16a L2}
\cos\bar\Phi\, |J^{i,j}| = \frac{2}{km}\ (1 + o(1)),
 \eeq
 \beq\lab{16b L2}
\cos\bar\Phi\marka\, |J^{\s(i),j}| = \frac{2}{km}\ (1 + o(1)),
 \eeq
with $o(1)$ being uniformly small over all $i$,\, $j$ as $\ve \to
0$, and $|J|$ being the Lebesgue measure of $J$. (Recall that the
parameters $\bar\Phi$,\, $\bar\Phi\marka$,\, $k$ and the intervals
$J^{i,j}$ implicitly depend on $\ve$.)

Choose closed intervals $\tilde{J}^{i,j} \subset {J}^{i,j}$ and
$\tilde{J}^{\s(i),j} \subset {J}^{\s(i),j}$ satisfying the following
conditions:~ (i) $\vphi\marka(\tilde{J}^{i,j} \times \{ 0 \}) =
\tilde{J}^{\s(i),j};\,$ (ii) some neighborhoods of $\tilde{J}^{i,j}$
and $\tilde{J}^{\s(i),j}$ belong to ${J}^{i,j}$ and ${J}^{\s(i),j}$,
respectively;\, and (iii) the pseudo-billiard trajectory with the
initial data $(\vphi, 0)$,\, $\vphi \in \tilde{J}^{i,j}$ does not
intersect the neighbor lateral reflectors (that is, the lateral
reflectors corresponding to the intervals $\tilde{J}^{i,j+1}$ and
$\tilde{J}^{\s(i),j-1}$, if $j \ne 1,\ k$; if $j = 1$ then
$\tilde{J}^{\s(i),j-1}$ should be replaced with
$\tilde{J}^{\s(i)-1,k}$, and if $j = k$ then $\tilde{J}^{i,j+1}$
should be replaced with $\tilde{J}^{i+1,1}$). Note in this regard
that the neighbor lateral reflectors occupy a small part of the
angular intervals ${J}^{i,j}$ and ${J}^{\s(i),j}$ (represented on
the figure by the arcs $AB$ and $B\marka A\marka$). Other lateral
reflectors will not be intersected, by the choice of lateral
reflectors.

By virtue of (\ref{15 L2}),\, (\ref{16a L2}),\, (\ref{16b L2}) and
because of smallness of the angular intervals occupied by the
lateral reflectors, $\tilde{J}^{i,j}$ and $\tilde{J}^{\s(i),j}$ may
be chosen in such a way that the ratios $\lam(\tilde{J}^{i,j}) /
\lam({J}^{i,j})$ and $\lam(\tilde{J}^{\s(i),j}) /
\lam({J}^{\s(i),j})$ uniformly (with respect to $i$,\, $j$) tend to
1 as $\ve \to 0$. Thus, a billiard particle going from $O$ in a
direction $\vphi \in \tilde{J}^{i,j}$, makes the same reflections
and in the same order as under the pseudo-billiard dynamics: first,
reflection from $\EEE_{i,j}$, then from $\PPP_{i,j}$, from
$\PPP\marka_{i,j}$, from $\EEE\marka_{i,j}$, and finally, the
particle goes back to $O$ in the direction $\vphi\marka(\vphi, 0)
\in \tilde{J}^{\s(i),j}$.

Choose $a_{i,j}$ in such a way that the following conditions are
fulfilled: if $(\xi, \vphi) \in [-a_{i,j},\, a_{i,j}] \times \tilde
J^{i,j}$, then~ (i) the corresponding billiard trajectory does not
intersect the lateral reflectors and the indicated order of
reflections is preserved;~ (ii) $\vphi\marka(\vphi, \xi) \in
J^{\s(i),j}$;~ (iii) $\left| \frac{\pl\xi\marka}{\pl\xi}(\vphi, \xi)
\right| \in [(1 + \ve)^{-3},\ (1 + \ve)^3]$. Analogously, choose
$a_{\s(i),j}$ in such a way that the conditions are fulfilled: if
$(\xi, \vphi) \in [-a_{\s(i),j},\, a_{\s(i),j}] \times \tilde
J^{\s(i),j}$ then~ (i) he billiard trajectory does not intersect the
lateral reflectors and the order of its reflections is reverse;~
(ii) $\vphi\marka(\vphi, \xi) \in J^{i,j}$;~ (iii) $\left|
\frac{\pl\xi\marka}{\pl\xi}(\vphi, \xi) \right| \in [(1 +
\ve)^{-3},\ (1 + \ve)^3]$. Note that the values $a_{i,j}$ and
$a_{\s(i),j}$ implicitly depend on $\ve$.

Select $a_\ve \le \min_{i,j} a_{ij}$ in such a way that $a_\ve \to
0$ as $\ve \to 0$ and denote $I_\ve = (-a_\ve,\, a_\ve) \times \{ 0
\}$,\, $\tilde I_\ve = (-a_\ve(1 + \ve)^{-3},\, a_\ve(1 + \ve)^{-3})
\times \{ 0 \}$, and $\tilde J^i_\ve = \tilde J^i := \cup_j \tilde
J^{i,j}_\ve$. The part of the boundary of $\Om_\ve$ related to the
angular intervals ${J}^{i,j}$ and ${J}^{\s(i),j}$ under
consideration is formed by the arcs of ellipses $\EEE_{i,j}$,\,
$\EEE\marka_{\!\!i,j}$ and the corresponding lateral reflectors.
Then a billiard particle with initial conditions $(\xi, \vphi) \in
\tilde I_\ve \times \tilde J^{i,j}$ after making four reflections
will intersect $l$ at a point $(\xi\marka, 0) \in I_\ve$, and the
angle at the point of intersection will be $\vphi^+_{\Om_\ve,I_\ve}
(\xi, \vphi) = \vphi\marka(\vphi, \xi) \in {J}^{\s(i),j} \subset
J^{\s(i)}$. Thus, the set of values $(\xi, \vphi) \in I_\ve \times
J^i$ such that $\vphi^+_{\Om_\ve,I_\ve} (\xi, \vphi) \not\in
J^{\s(i)}$ is contained in the set $\left( I_\ve \times J^i \right)
\setminus \left( \tilde I_\ve \times \tilde J^i_\ve \right)$, whose
measure is vanishing as $\ve \to 0$. Q.E.D.

  {\bf (b)}~  If $2 \le i = \s(i) \le m - 1$, the corresponding part of the boundary
is the arc of circumference of radius 2 with the center at $O$
occupying the angular interval $J^i$, that is, the set $\{
2(\sin\vphi, \cos\vphi),\ \vphi \in J^i \}$. Next we will show that
for all values $(\xi, 0) \in I_\ve$,\, $\vphi \in J^i$, except for a
portion of order $o(1)$, the corresponding billiard particle makes
one reflection from the arc and then goes back to $I_\ve$ in the
direction $\vphi\marka \in J^i$.

For all values $\vphi \in J^i$, except for the union of two
intervals of vanishing length (each of the intervals is contained in
$J^i$, has the length $2\arctan(a_\ve/4)$, and contains an endpoint
of $J^i$), the particle starting at $(\xi, 0) \in I_\ve$ in the
direction $\vphi$ will reflect from the indicated arc of
circumference. Let $\psi \in J^i$ be the angular coordinate of the
reflection point. By $(\xi\marka, 0)$ denote the point at which the
reflected particle intersects the straight line $l$. One easily
verifies that
 \beq\lab{pssi}
  \frac{1}{\xi}\ +\ \frac{1}{\xi\marka}\ =\ \cos\psi.
 \eeq
One has
  \beq\lab{kssi}
|\xi|\ <\ a_\ve,
 \eeq
hence
 \beq\lab{pssipssi}
\frac{1}{|\xi\marka|}\ =\ \big| \cos\psi - \frac{1}{\xi}\, \Big|\ >\
\frac{1}{a_\ve} - 1.
 \eeq
From (\ref{pssi}) it follows that ${|\xi + \xi\marka|}/{|\xi
\xi\marka|} = |\cos\psi| \le 1$, and taking into account
(\ref{kssi}) and (\ref{pssipssi}), one finds that $|\xi + \xi\marka|
< {a_\ve^2}/{(1 - a_\ve)}$. This implies that for all values $(\xi,
0) \in I_\ve$, except for a set of measure $O(a_\ve^2)$, the second
point of intersection of the billiard trajectory belongs to $I_\ve$,
moreover the velocity at this point, $\vphi^+_{\Om_\ve,I_\ve} (\xi,
\vphi)$, belongs to $\NNN_{2\arctan(a_\ve/4)}(J^i)$, the
neighborhood of $J^i$ of radius $2\arctan(a_\ve/4)$. This finally
implies that for all $(\xi, \vphi) \in I_\ve \times J^i$, except for
a portion of order $O(a_\ve)$, holds $\vphi^+_{\Om_\ve,I_\ve} (\xi,
\vphi) \in J^i$.

{\bf (c)}~ The parts of the hollow's boundary, corresponding to
$J^1$ and $J^m$, are formed by smooth curves joining the
corresponding endpoints of $I_\ve$ and the points $2(\sin\Phi_0,
-\cos\Phi_0)$ and $2(\sin\Phi_0, \cos\Phi_0)$, respectively. The
unique condition on these curves is that they can be parametrized by
the monotonically increasing angular coordinate. For those values
$\s(1)$,\, $\s(m)$ that do not coincide with neither 1 nor $m$ take
just the arcs of circumference of radius 2 corresponding to the
angular intervals $J^{\s(1)}$,\, $J^{\s(m)}$.

Consider the union of all the elliptic arcs $\EEE_{i,j}$,\,
$\EEE\marka_{\!\!i,j}$ introduced in item (a), all the arcs of
circumference defined in items (a) and (b), and the two curves
introduced in this item (c). Let us call this union the {\it main
element}. Each lateral reflector is a curve; select it in such a way
that both its endpoints belong to the main element. Finally, the
curve $\pl\Om_\ve \setminus I_\ve$ is the union of all the lateral
reflectors and the part of the main element visible from $O$ (that
is, which is not shielded by the adjacent lateral reflectors). Thus,
the definition of the hollow $(\Om_\ve, I_\ve)$ is complete.

On the figure below, there is shown a particular hollow $(\Om_\ve,
I_\ve)$ corresponding to the permutation $\s = \big(\! {\scriptsize
\begin{array}{ccccc} 1 & 2 & 3 & 4 & 5\\ 5 & 4 & 3 & 2 & 1
\end{array}}\! \big)$. The angular intervals $J^1, \ldots, J^5$ are
separated by dotted lines. The family of hollows $(\Om_\ve, I_\ve)$,
with vanishingly small $\ve$, has the following property: for almost
all particles with the initial direction from $J^2$ (resp. $J^3$,\,
$J^4$), the final direction will belong to $J^4$ (resp. $J^3$,\,
$J^2$). On the figure, there is shown the trajectory of a particle
with the initial direction $\vphi \in J^2$ and the final direction
$\vphi^+ \in J^4$. The particle makes a reflection from an elliptic
arc, then two reflections from (very small) parabolic arcs, and
finally, again from an elliptic arc. According to our notation,
these arcs are $\EEE_{2,2}$,\, $\PPP_{2,2}$,\,
$\PPP\marka_{\!\!2,2}$, and $\EEE\marka_{\!\!2,2}$.
 \newpage

 
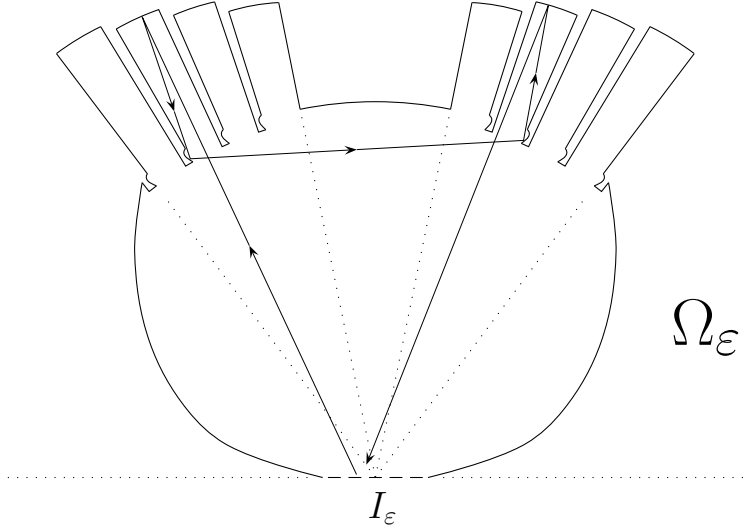
\vspace*{66mm}

 \rput(6,1){
 \psarc[linewidth=0.4pt](0,0){5}{78.46}{101.54}
 \psline[linestyle=dashed,linewidth=0.4pt](-0.63,0)(0.63,0)
  \psline[linestyle=dotted,linewidth=0.5pt](-0.9,0)(-5,0)
  \psline[linestyle=dotted,linewidth=0.5pt](0.9,0)(5,0)
  \psline[linestyle=dotted,linewidth=0.5pt](0,0)(-1,4.899)
  \psline[linestyle=dotted,linewidth=0.5pt](0,0)(-2.82,3.76)
  \psline[linestyle=dotted,linewidth=0.5pt](0,0)(1,4.899)
  \psline[linestyle=dotted,linewidth=0.5pt](0,0)(2.82,3.76)
 \pscurve[linewidth=0.4pt](0.7,0)(2,0.45)(2.5,0.9)(3,1.8)(3.2,3)(3.1,3.923)
 \pscurve[linewidth=0.4pt](-0.7,0)(-2,0.45)(-2.5,0.9)(-3,1.8)(-3.2,3)(-3.1,3.923)
    \psline[linewidth=0.4pt](3.1,3.923)(3.007,3.805)(2.91,3.88)
    \psline[linewidth=0.4pt]       (3.03,4.04)(4.24,5.64)
    \psline[linewidth=0.4pt](3.653,6)(2.522,4.1427)(2.425,4.2)
    \psline[linewidth=0.4pt](2.525,4.3733)    (3.445,5.967)
    \psline[linewidth=0.4pt](2.882,6.228)(2.037,4.4015)(1.94,4.445)
    \psline[linewidth=0.4pt]    (2.02,4.6284)(2.684,6.15)
    \psline[linewidth=0.4pt](2.137,6.327)(1.552,4.595)(1.455,4.627)
    \psline[linewidth=0.4pt]    (1.515,4.8177)(1.95,6.2)
    \psline[linewidth=0.4pt](1.285,6.321)(1,4.899)
    \psline[linewidth=0.4pt](-3.1,3.923)(-3.007,3.805)(-2.91,3.88)
    \psline[linewidth=0.4pt]       (-3.03,4.04)(-4.24,5.64)
    \psline[linewidth=0.4pt](-3.653,6)(-2.522,4.1427)(-2.425,4.2)
    \psline[linewidth=0.4pt](-2.525,4.3733)    (-3.445,5.967)
    \psline[linewidth=0.4pt](-2.882,6.228)(-2.037,4.4015)(-1.94,4.445)
    \psline[linewidth=0.4pt]    (-2.02,4.6284)(-2.684,6.15)
    \psline[linewidth=0.4pt](-2.137,6.327)(-1.552,4.595)(-1.455,4.627)
    \psline[linewidth=0.4pt]    (-1.515,4.8177)(-1.95,6.2)
    \psline[linewidth=0.4pt](-1.285,6.321)(-1,4.899)
  \psarc[linewidth=0.4pt](2.12,2.82){3.528}{53.13}{64.25}
  \psarc[linewidth=0.4pt](1.814,3.142){3.262}{60}{71}
  \psarc[linewidth=0.4pt](1.49,3.415){2.984}{66.42}{77.5}
  \psarc[linewidth=0.4pt](1.154,3.67){2.652}{72.54}{87.1}
  \psarc[linewidth=0.4pt](-2.12,2.82){3.528}{115.75}{126.87}
  \psarc[linewidth=0.4pt](-1.814,3.142){3.262}{109}{120}
  \psarc[linewidth=0.4pt](-1.49,3.415){2.984}{102.5}{113.58}
  \psarc[linewidth=0.4pt](-1.154,3.67){2.652}{92.9}{107.46}
     \psecurve[linewidth=0.4pt](2.975,4.05)(3.03,4.04)(3.04,3.97)(2.91,3.88)(2.79,3.845)
     \psecurve[linewidth=0.4pt](-2.975,4.05)(-3.03,4.04)(-3.04,3.97)(-2.91,3.88)(-2.79,3.845)
     \psecurve[linewidth=0.4pt](2.4692,4.3767)(2.525,4.3733)(2.5433,4.305)(2.425,4.2)(2.31,4.151)
     \psecurve[linewidth=0.4pt](-2.4692,4.3767)(-2.525,4.3733)(-2.5433,4.305)(-2.425,4.2)(-2.31,4.151)
     \psecurve[linewidth=0.4pt](1.9641,4.6256)(2.02,4.6284)(2.0459,4.5626)(1.94,4.445)(1.8311,4.3831)
     \psecurve[linewidth=0.4pt](-1.9641,4.6256)(-2.02,4.6284)(-2.0459,4.5626)(-1.94,4.445)(-1.8311,4.3831)
     \psecurve[linewidth=0.4pt](1.46,4.809)(1.515,4.8177)(1.5477,4.755)(1.455,4.627)(1.3535,4.554)
     \psecurve[linewidth=0.4pt](-1.46,4.809)(-1.515,4.8177)(-1.5477,4.755)(-1.455,4.627)(-1.3535,4.554)
   \psline[linewidth=0.3pt,arrows=->,arrowscale=1.6](-0.25,0.04)(-1.675,3.085)
   \psline[linewidth=0.3pt,arrows=->,arrowscale=1.6](-1.675,3.085)(-3.1,6.13)(-2.6733,4.87)
   \psline[linewidth=0.3pt,arrows=->,arrowscale=1.6](-2.6733,4.87)(-2.46,4.24)(-0.245,4.3625)
   \psline[linewidth=0.3pt,arrows=->,arrowscale=1.6](-0.245,4.3625)(1.97,4.485)(2.135,5.3825)
   \psline[linewidth=0.3pt,arrows=->,arrowscale=1.6](2.135,5.3825)(2.3,6.28)(-0.12,0.17)
   \rput(0.1,-0.4){\large $I_\ve$}
   \rput(4.4,2){\Huge $\Om_\ve$}
 }
  \vspace{0mm}

\section{Concluding remarks and applications}

Physical bodies in the real world have atomic structure and
therefore are disconnected. This is a reason for using (generally)
disconnected sets $Q_m$ in the definition of a rough body. In future
we intend to turn to propose and study the notion of a
three-dimensional rough body, where the connectivity assumption is
absolutely useless; this is another reason. By removing this
assumption, the consideration in two dimensions (namely, proof of
lemma \ref{l1}) is made somewhat more difficult, but at the same
time prerequisites for passing to the three-dimensional case are
created.

In fact, the notions of "disconnected"{} (as everywhere in this
paper) and "connected"{} rough bodies are equivalent. There is a
natural one-to-one correspondence between the equivalence classes in
the connected and disconnected cases,\footnote{more precisely,
equivalence classes formed by sequences of
connected\,/\,disconnected sets} the former classes being subclasses
of the latter ones under this correspondence.

Let us now consider applications of theorem to problems of the body
of minimal or maximal aerodynamic resistance. A two-dimensional
convex body $B$ moves, at constant velocity, through a rarefied
homogeneous medium in $\RRR^2$, and at the same time slowly rotates.
The rotation is generally non-uniform; we assume that during a
sufficiently long observation period, in a reference system
connected with the body the body's velocity is distributed in $S^1$
according to a given density function $\rho$, with $\int_{S^1}
\rho(v)\, dv = 1$. The medium particles do not mutually interact,
and collisions of the particles with the body are absolutely
elastic. The resistance of the medium to the motion of the body is a
vector-valued function of time. After averaging it over a
sufficiently long period of time, one gets a vector. We are
interested in the projection of this vector on the direction of
motion; for the sake of brevity, it will be called mean resistance,
or just resistance. The problem is: given $B$, determine the
roughness on it in such a way that main resistance of the resulting
rough body is minimal or maximal.

A prototype of such a mechanical system is an artificial satellite
of the Earth on relatively low altitudes ($100 \div 200$ km), with
restricted capacity of rotation angle control. The satellite's
motion is slowing down by the rest of atmosphere; the problem is
minimize or maximize the effect of slowing down. The problems of
resistance {\it maximization} may also arise when considering solar
sail: a spacecraft driven by the pressure of solar photons.

The initial velocity of an incident particle (in the reference
system connected with the body) is $-v$, and the final velocity is
$v^+$; therefore, the momentum transmitted by the particle to the
body is $v + v^+$. The projection of the transmitted momentum on the
direction of motion of the body equals $1 + \langle v,\, v^+
\rangle$. Averaging this value over all particles incident on the
body within a sufficiently long time interval, one gets the mean
resistance. The averaging amounts to integration over $\rho(v)\,
d\nu_\BBB(v, v^+, n)$; that is, mean resistance of the rough body
equals
$$
R(\nu_\BBB)\ =\ \int\!\!\!\int\!\!\!\int_{\TTT^3} (1 + \langle v,\,
v^+ \rangle)\, \rho(v)\, d\nu_\BBB(v, v^+, n).
$$
Using theorem 1 and Fubini's theorem, one rewrites this formula in
the form
 \beq\lab{resis}
R(\nu_\BBB)\ =\ \int_{S^1} d\tau_B(n) \int\!\!\!\int_{\TTT^2} (1 +
\langle v,\, v^+ \rangle)\, \rho(v)\, d\eta_{\BBB,n}(v, v^+),
 \eeq
where $\eta_{\BBB,n} \in \Lam_n$. Thus, the minimization problem for
$R(\nu_\BBB)$ reduces to minimization, for any $n$, of the
functional $\int\!\!\!\int_{\TTT^2} (1 + \langle v,\, v^+ \rangle)\,
\rho(v)\, d\eta(v, v^+)$ over all $\eta \in \Lam_n$. Using the
notation introduced in subsection 3.3, one comes to the problem:
 \beq\lab{reduced problem}
\inf_{\eta \in \Lam} \int\!\!\!\!\int_{\Box} (1 + \cos(\vphi -
\vphi^+))\, \varrho(\vphi) \, d\eta(\vphi, \vphi^+),
 \eeq
where $\varrho(\vphi) = \rho(v)$ for $\vphi = \text{Arg}\,v -
\text{Arg}\,n$. This problem, in turn, by symmetrization of the cost
function reduces to a particular Monge-Kantorovich problem:
 \beq\lab{MK problem}
\inf_{\eta \in \Lam_{\lam,\lam}} \FFF(\eta),~~~ \text{where}~~~
\FFF(\eta)\ =\ \int\!\!\!\!\int_{\Box} c(\vphi,\vphi^+)\,
d\eta(\vphi, \vphi^+),
 \eeq
where $c(\vphi, \vphi^+) = (1 + \cos(\vphi - \vphi^+))\,
\frac{\varrho(\vphi) + \varrho(\vphi^+)}{2}$ and $\Lam_{\lam,\lam}$
is the set of measures $\eta$ on $\Box$ having both marginal
measures equal to $\lam$:\, \ $\pi_\vphi^\# \eta = \lam =
\pi_{\vphi^+}^\# \eta$. Recall that $\lam$ is defined by $d\lam(v) =
\cos\vphi\, d\vphi$.

The problem (\ref{MK problem}) can be exactly solved in several
particular cases. Consider the case of uniform motion, where the
function $\rho$, and therefore $\varrho$, is constant, and thus, one
can take $c(\vphi,\vphi^+) = \frac 38\, (1 + \cos(\vphi -
\vphi^+))$.\footnote{The normalization constant $3/8$ is taken for
further convenience.} Note that $\FFF(\eta_0) = 1$ and therefore
resistance of the smooth body is equal to its perimeter: $R(\nu_B) =
\int_{S^1} d\tau_B(n)\, \FFF(\eta_0) = |\pl B|$. (Recall that the
measure $\eta_0$ belongs to $\Lam$ and is supported on the diagonal
$\vphi^+ = -\vphi$.) The minimization problem (\ref{MK problem}) for
constant $\varrho$ was solved in \cite{sb-math-averaged04}: one has
$\inf_{\BBB} R(\nu_\BBB) = 0.9878... \cdot |\pl B|$, the infimum
being taken over all roughenings of $B$.

Note that the corresponding {\it maximization} problem for (\ref{MK
problem}) has the trivial solution, which does not depend on the
function $\varrho$:\, $\eta = \eta_\star$, the measure $\eta_\star
\in \Lam$ being supported on the diagonal $\vphi^+ = \vphi$. One has
$\sup_{\BBB} R(\nu_\BBB) = \kappa\, |\pl B|$, where $\kappa = \left(
\int_{-\pi/2}^{\pi/2} \varrho(\vphi)\, \cos\vphi\, d\vphi \right)
\Big/ \left( \int_{-\pi/2}^{\pi/2} \varrho(\vphi)\, \cos^3 \vphi\,
d\vphi \right) > 1$; in the case of uniform rotation one has $\kappa
= 1.5$. The maximization problem was studied in more detail in
\cite{maximization}.

\section*{Appendix A}

The construction is simple (see the figure), but its description is
a bit cumbersome.

Take a point in the interior of $B$ and connect it by segments with
all vertices. The polygon is thus divided into several triangles;
fix $i$ and $m$ and consider the triangle with the base
$\mathfrak{b}_i$, the $i$th side of $B$. Denote by $d(\Om_i^m)$ the
diameter of the orthogonal projection of $\Om_i^m$ on the straight
line containing $I_i^m$; one obviously has $d(\Om_i^m) \ge |I_i^m|$.
Fix a positive number $\kappa < |I_i^m| / d(\Om_i^m)$. Take a
rectangle $\Pi^1$ contained in the triangle and such that one side
of $\Pi^1$ belongs to $\mathfrak{b}_i$. By $\del_1$ denote the total
length of the part of $\mathfrak{b}_i$ which is not occupied by this
side.

For the sake of brevity, the image of a set under the composition of
a homothety with positive ratio and a translation will be called a
copy of this set. Take several copies of $\Om_i^m$ (copies of first
order) that do not mutually interact, belong to $\Pi^1$, the
corresponding copies of $I_i^m$ belong to $\mathfrak{b}_i$, and the
portion of the side of $\Pi^1$ occupied by them is more than
$\kappa$.
 
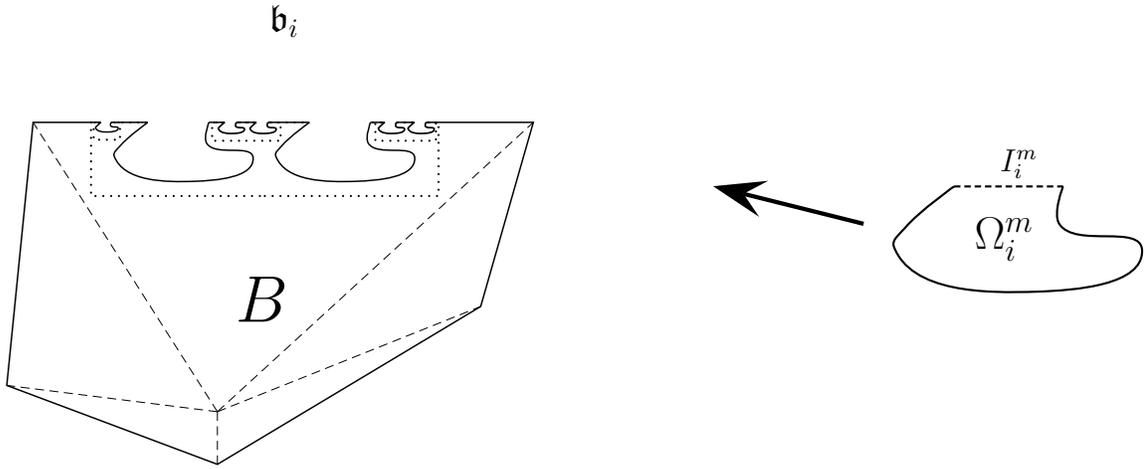
\vspace{59mm}

 \rput(-1,0){
 \scalebox{0.7}{
 \rput(4,0){
 \pspolygon(0,-1)(5,2)(6,5.5)(-3.5,5.5)(-4,0.5)
 \psline[linestyle=dashed,linewidth=0.4pt](0,0)(0,-1)
 \psline[linestyle=dashed,linewidth=0.4pt](0,0)(5,2)
 \psline[linestyle=dashed,linewidth=0.4pt](0,0)(6,5.5)
 \psline[linestyle=dashed,linewidth=0.4pt](0,0)(-3.5,5.5)
 \psline[linestyle=dashed,linewidth=0.4pt](0,0)(-4,0.5)
 \pspolygon[linestyle=dotted,linewidth=1.2pt](4.2,4.1)(4.2,5.5)(-2.4,5.5)(-2.4,4.1)
 \pspolygon[linestyle=dotted,linewidth=1.2pt](-1.85,5.17)(-1.85,5.5)(-2.35,5.5)(-2.35,5.17)
 \pspolygon[linestyle=dotted,linewidth=1.2pt](1.2,5.15)(1.2,5.5)(-0.1,5.5)(-0.1,5.15)
 \pspolygon[linestyle=dotted,linewidth=1.2pt](4.15,5.15)(4.15,5.5)(3,5.5)(3,5.15)
 }
 \rput(2.8,5.5){
 \scalebox{0.4}{
 \psline[linecolor=white,linewidth=4.8pt](-0.2,0)(2.7,0)
 \pscurve[linewidth=2pt](2.7,0)(2.6,-1)(4.8,-1.7)(2,-2.8)(-1.8,-1.6)(-1.5,-1.1)(-0.2,0)
 }
 }
 \rput(5.8,5.5){
 \scalebox{0.4}{
 \psline[linecolor=white,linewidth=4.8pt](-0.2,0)(2.7,0)
 \pscurve[linewidth=2pt](2.7,0)(2.6,-1)(4.8,-1.7)(2,-2.8)(-1.8,-1.6)(-1.5,-1.1)(-0.2,0)
 }}
 \rput(1.8,5.5){
 \scalebox{0.065}{
 \psline[linecolor=white,linewidth=40pt](-0.2,0)(2.7,0)
 \pscurve[linewidth=12pt](2.7,0)(2.6,-1)(4.8,-1.7)(2,-2.8)(-1.8,-1.6)(-1.5,-1.1)(-0.2,0)
 }}
 \rput(4.75,5.5){
 \scalebox{0.075}{
 \psline[linecolor=white,linewidth=40pt](-0.2,0)(2.7,0)
 \pscurve[linewidth=12pt](2.7,0)(2.6,-1)(4.8,-1.7)(2,-2.8)(-1.8,-1.6)(-1.5,-1.1)(-0.2,0)
 }}
 \rput(4.15,5.5){
 \scalebox{0.075}{
 \psline[linecolor=white,linewidth=40pt](-0.2,0)(2.7,0)
 \pscurve[linewidth=12pt](2.7,0)(2.6,-1)(4.8,-1.7)(2,-2.8)(-1.8,-1.6)(-1.5,-1.1)(-0.2,0)
 }}
 \rput(7.75,5.5){
 \scalebox{0.075}{
 \psline[linecolor=white,linewidth=40pt](-0.2,0)(2.7,0)
 \pscurve[linewidth=12pt](2.7,0)(2.6,-1)(4.8,-1.7)(2,-2.8)(-1.8,-1.6)(-1.5,-1.1)(-0.2,0)
 }}
 \rput(7.19,5.5){
 \scalebox{0.075}{
 \psline[linecolor=white,linewidth=40pt](-0.2,0)(2.7,0)
 \pscurve[linewidth=12pt](2.7,0)(2.6,-1)(4.8,-1.7)(2,-2.8)(-1.8,-1.6)(-1.5,-1.1)(-0.2,0)
 }}
 }
 \rput(3.3,1.5){\Huge $B$}
 \rput(3.6,5.2){\large $\mathfrak{b}_i$}
  \rput(-0.5,0){
  \rput(13,3){
    \scalebox{0.5}{
 \psline[linestyle=dashed,linewidth=1.6pt](-0.2,0)(2.7,0)
 \pscurve[linewidth=1.6pt](2.7,0)(2.6,-1)(4.8,-1.7)(2,-2.8)(-1.8,-1.6)(-1.5,-1.1)(-0.2,0)
 }
 }
 \rput(13.6,2.3){\Large $\Om_i^m$}
 \rput(13.8,3.3){$I_i^m$}
 }
  \psline[linewidth=1.6pt,arrows=<-,arrowscale=3](9.3,3)(11.3,2.5)
  }
  \vspace{15mm}

Next, take several rectangles that do not mutually intersect and do
not intersect with the chosen copies of $\Om_i^m$, belong to
$\Pi^1$, and have one side contained in $\mathfrak{b}_i$. Denote by
$\Pi^2$ the union of these rectangles and by $\del_2$, the total
length of the part of the side of $\Pi^1$ which is not occupied by
the rectangles from $\Pi^2$ and by the copies of $I_i^m$. Next, for
each rectangle from $\Pi^2$ choose several copies of $\Om_i^m$
(copies of second order) in the way completely similar to the
described above (see fig.\,\ref{fig:apa}).

Continuing this process, one obtains a sequence $\Pi^1$,\, $\Pi^2,
\ldots$ of unions of rectangles and collections of copies of
$\Om_i^m$ of 1st, 2nd, $\ldots$ order. Choose the rectangles in such
a way that $\del_1 + \del_2 + \ldots < 1/m$ and Area$(\Pi^1) < 1/m$.
Finally, choose $k$ such that the total length of sides of
rectangles from $\Pi^{k+1}$ contained in $\mathfrak{b}_i$ is less
than $1/m$, and take the collection of copies of $\Om_i^m$ of order
$1,\, 2, \ldots, k$ (we shall call it full collection). The total
length of the part of $\mathfrak{b}_i$ not occupied by the
corresponding copies of $I_i^m$ is less than $2/m$, and therefore
goes to zero as $m \to \infty$.

By definition, the desired set $Q_m$ is $B$ minus the union of full
collections of copies of $\Om_i^m$ over all $i$.

\section*{Appendix B}

We prove here slightly more than needed.

\begin{predl}\label{u3}
Let $A = (a_{ij})_{i,j=1}^k$ be a symmetric matrix, with $a_{ij}$
being nonnegative integers. Denote $n_i = \sum_{j=1}^k a_{ij}$. Then
there exist matrices $B_{ij} = (b_{ij}^{\mu\nu})_{\mu,\nu}$ of size
$n_i \times n_j$ such that $b_{ij}^{\mu\nu} \in \{ 0,\, 1 \}$,\,
$B_{ij}^T = B_{ji}$, the sum of elements in $B_{ij}$ equals
$a_{ij}$, and the block matrix $D = (B_{ij})$ contains exactly one
unit in each row and each column.
\end{predl}

Note that for some values $i = i_1,\, i_2, \ldots$ it may happen
that $n_i = 0$, that is, $a_{ij} = 0$ for all $j = 1, \ldots, k$.
Then the corresponding matrices $B_{ij}$ have the size $0 \times
n_j$, that is, are empty. In this case $D$ coincides with the block
matrix $D' = (B_{ij})$ having the rows $i_1,\ i_2, \ldots$ and
columns $i_1,\ i_2, \ldots$ crossed out.

\begin{proof}
The proof is by induction on $k$. Let the statement be true for $k -
1$; prove it for $k$. Take the matrix $\tilde A = (a_{ij})_{i,j=
2}^k$; there exists a block matrix $\tilde B = (\tilde
B_{ij})_{i,j=2}^k$ satisfying the statement. Note that the order of
$\tilde B_{ij}$ is $\tilde n_i \times \tilde n_j$, where $\tilde n_i
= \sum_{j=2}^k a_{ij} = n_i - a_{i1}$. Define the matrices $B_{ij}$
as follows.

(a) Put $B_{11} =$ diag\,$\{ \underbrace{1, \ldots, 1}_{a_{11}},\ 0,
\ldots, 0 \}$.

(b) Put $b_{12}^{a_{11}+1,1} = \ldots =
b_{12}^{a_{11}+a_{12},a_{12}} = 1$;\ $b_{13}^{a_{11}+a_{12}+1,1} =
\ldots = b_{13}^{a_{11}+a_{12}+a_{13},a_{13}} = 1$;\ \ldots; \
$b_{1k}^{a_{11}+\ldots+a_{1,k-1}+1,1} = \ldots =
b_{1k}^{a_{11}+\ldots+a_{1k},a_{1k}} = 1$;\ the other elements of
the matrices $B_{1j}$,\, $j = 2, \ldots, k$ are zeros. Thus, on the
diagonal of $B_{1j}$ starting from the element at the first column
and the $(a_{11} + a_{12} + \ldots + a_{1,j-1} + 1)$th row, the
first $a_{1j}$ elements equal 1, and the remaining elements on this
diagonal and all the elements off the diagonal are zeros. This
defines the matrices $B_{1j}$,\, $j = 2, \ldots, k$. The matrices
$B_{i1}$,\, $i = 2, \ldots, k$ are determined by the condition
$B_{i1} = B_{1i}^T$.

(c) For $i \ge 2$,\, $j \ge 2$ define the matrix $B_{ij}$ as
follows. For $\mu \le a_{1i}$ or $\nu \le a_{1j}$, put
$b_{ij}^{\mu\nu} = 0$, and for $\mu \ge a_{1i} + 1$,\, $\nu \ge
a_{1j} + 1$, put $b_{ij}^{\mu\nu} = \tilde
b_{ij}^{\mu-a_{1i},\nu-a_{1j}}$. Thus, in the obtained matrix
$B_{ij}$, the right lower corner coincides with the matrix $\tilde
B_{ij}$, and all the remaining elements are equal to zero. The
number of rows of this matrix equals $a_{1i} + \tilde n_i = n_i$,
and the number of columns equals $a_{1j} + \tilde n_j = n_j$. One
obviously has $B_{ij}^T = B_{ji}$.

One easily verifies that $\sum_{\mu\nu} b_{ij}^{\mu\nu} = a_{ij}$
and that each row and each column of the obtained block matrix $D =
(B_{ij})_{i,j=1}^k$ contains precisely one unit.

\end{proof}

\section*{Acknowledgements}

This work was supported by {\it Centre for Research on Optimization
and Control} (CEOC) from the ''{\it Funda\c{c}\~{a}o para a
Ci\^{e}ncia e a Tecnologia}'' (FCT), cofinanced by the European
Community Fund FEDER/POCTI, and by FCT (research project
PTDC/MAT/72840/2006).

\end{document}